\documentclass[12pt,reqno]{amsart}
\usepackage{amsmath}
\usepackage{amssymb}
\usepackage{amsthm}
\usepackage{mathrsfs}
\usepackage{latexsym}

\newtheorem{theorem}{Theorem}[section]
\newtheorem{lemma}{Lemma}[section]
\newtheorem{corollary}{Corollary}[section]
\newtheorem{remark}{Remark}[section]
\newtheorem{definition}{Definition}[section]
\newtheorem{proposition}{Proposition}[section]
\newtheorem{example}{Example}[section]
\newtheorem{assumption}{Assumption}[section]
\numberwithin{equation}{section}
\newcommand{\bth}{\begin{theorem}}
\newcommand{\ethe}{\end{theorem}}

\newcommand{\bre}{\begin{remark}}
\newcommand{\ere}{\end{remark}}

\newcommand{\ble}{\begin{lemma}}
\newcommand{\ele}{\end{lemma}}
\newcommand{\bde}{\begin{definition}}
\newcommand{\ede}{\end{definition}}

\newcommand{\bco}{\begin{corollary}}
\newcommand{\eco}{\end{corollary}}

\newcommand{\bpr}{\begin{proposition}}
\newcommand{\epr}{\end{proposition}}

\newcommand{\bexer}{\begin{exercise}}
\newcommand{\eexer}{\end{exercise}}
\newcommand{\breh}{\begin{hint}}
\newcommand{\ereh}{\end{hint}}

\newcommand{\halmos}{\hfill \qed}

\newcommand{\bexam}{\begin{example}}
\newcommand{\eexam}{\end{example}}

\newcommand{\pr} {{\bf Proof.}}

\newcommand{\bfi}{\begin{fig}}
\newcommand{\efi}{\end{fig}}

\newcommand{\beao}{\begin{eqnarray*}}
\newcommand{\eeao}{\end{eqnarray*}\noindent}

\newcommand{\beam}{\begin{eqnarray}}
\newcommand{\eeam}{\end{eqnarray}\noindent}
\newcommand{\E}{\mathbf{E}}
\newcommand{\PP}{\mathbf{P}}
\newcommand{\nto}{n\to\infty}

\newcommand{\xto}{x\to\infty}
\newcommand{\tto}{t\to\infty}

\newcommand{\bF}{\overline{F}}

\newcommand{\bV}{\overline{V}}

\newcommand{\bbr}{{\mathbb R}}

\newcommand{\bbb}{{\mathbb B}}

\newcommand{\bbn}{{\mathbb N}}

\newcommand{\vep}{\varepsilon}

\begin{document}
\title[Heavy-tailed random vectors: theory and applications]{Heavy-tailed random vectors: theory and applications}

\author[ D.G. Konstantinides, C. D. Passalidis ]{ Dimitrios G. Konstantinides, Charalampos  D. Passalidis}

\address{Dept. of Statistics and Actuarial-Financial Mathematics,
University of the Aegean,
Karlovassi, GR-83 200 Samos, Greece}
\email{konstant@aegean.gr,\;sasd24009@sas.aegean.gr.}

\date{{\small \today}}

\begin{abstract}
In this paper we introduce and study several multivariate heavy-tailed distribution classes, and we explore their closure properties and their applications. We consider the class of multivariate positively decreasing distributions, and its intersection with other multivariate distribution classes. Next, we show that the smallest of these new classes contains the standard multivariate regular variation class under a condition on Radon measure. We provide necessary and sufficient conditions for the closure property with respect to convolution in the class of multivariate subexponential positively decreasing  distributions. The conditional closure property with respect to convolution roots, in the class of multivariate subexponential distributions is also studied. Further, we study the closure properties with respect to scale mixtures, under the assumption that the random variable, that produces the scale mixture is weakly dependent to primary random vector. We give also a multivariate, dependent version of the Breiman's theorem. Under similar conditions, we establish the single big jump principle in scale mixture sums, under several conditions on the distributions of primary random vectors. Next, we present the asymptotic estimation of the finite time ruin probability, in a multivariate, time-dependent Poisson risk model. More precisely, we consider constant interest force, and common Poisson counting process, with any claim vector to be weakly dependent with inter-arrival times. Finally we give the lower bound of vector type precise large deviations. In case of non-random sums, the lower bound is estimated uniformly, the random vectors have some weak dependence structure, and we do not assume some condition for the distribution of the random vectors. In case of random sums, we provide similar estimations, under the additional condition that the random vectors are weak-equivalent on the chosen rare set. An application to compound risk model is also provided.
\end{abstract}

\maketitle
\vspace{3mm}



\section{Introduction} \label{sec.KP.1}

The heavy-tailed distributions, together with the dependence modeling occupy, without doubt, the most of the interest in the theoretical and applied probability. A flexible framework, in which these two topics are combined, represent the characterization of several classes of heavy-tailed distributions of random vectors. We find out that, while the characterization of such classes under terms of extreme value theory is well-established, see for example in 
\cite[Ch. 8]{beirlant:goegebeur:segers:teugels:2004}, this is not true for the multivariate distributions with heavy tails, except the multivariate 
regular variation, $MRV$, which is common in both approaches. As was mentioned in \cite{goldie:resnick:1988}, for the one-dimensional cases, beyond the regular variation, there are several differences between these two approaches. Similarly, these differences remain in multivariate distributions, for example if the distribution of a random vector ${\bf X}$ belongs to the class of multivariate Gumbel distributions and has asymptotically dependent components, then the marginal distributions of the random vector should belong to class of rapidly varying distributions, see \cite{asimit:furman:tang:vernic:2011}. That means, all the dominatedly varying distributions are excluded for the marginal ones, that contain important distributions, not restricted to the regular variation only. This way, we can understand that in some applications, as in risk theory, queuing theory but also in risk management, the multivariate extreme value theory does not suffice by itself, but it needs a complement study for multivariate heavy-tailed distributions.

A problem appearing in the study of multivariate heavy-tailed distributions, is the definition of multivariate subexponentiality. Although there are already at 
least four different such definitions, given by \cite{cline:resnick:1992}, \cite{omey:2006}, \cite{samorodnitsky:sun:2016} and \cite{konstantinides:passalidis:2024c}, with all of them to have found applications, however in the case of first three there is focus on the single big jump, while in the last one on the max-sum equivalence of the joint distribution tail, for the bivariate set up, see the discussion in \cite[Sec. 1]{konstantinides:passalidis:2024g}. 

In this paper we follow the approach by \cite{samorodnitsky:sun:2016}, with respect to the definition of multivariate subexponentiality and we have two targets. The first is to study several properties of multivariate heavy-tailed distributions, as the closure properties, the multivariate linear single big jump principle or the precise large deviations in multivariate set up. The 
second one is to present applications, as for example in risk theory, where these classes demonstrate a strong advantage in generalization of already existing risk models.

The paper is organized as follows. In Section 2, some preliminaries are provided for heavy-tailed distributions and next, following 
\cite{samorodnitsky:sun:2016} and \cite{konstantinides:passalidis:2024g}, is introduced the class of multivariate positively decreasing distributions and its intersection with other multivariate classes. In Section 3, we give some first results in relation with the ordering of multivariate classes and their closure properties, as the closure with respect to strong equivalence and with respect to convolution of random vectors, for classes related to multivariate positively decreasing distributions. We provide also a result of conditional closure property of the class of multivariate subexponential distributions with respect to convolution roots. In Section 4, we concentrate our attention 
on the scale mixtures and the scale mixture sums. As scale mixture we have in mind the quantity
\beam \label{eq.KP.1.1}
\Theta^{(i)}\,{\bf X}^{(i)} :=(\Theta^{(i)}\,X_{1}^{(i)},\,\ldots,\,\Theta^{(i)}\,X_{d}^{(i)})^{\top}\,,
\eeam
where with ${\bf x}^{\top}$ we denote the transpose of vector ${\bf x}$. The modeling through \eqref{eq.KP.1.1} find numerous practical applications, especially in actuarial mathematics and risk management, see Section 4.1 for details. Motivated by practical needs, we consider that in the scale mixture \eqref{eq.KP.1.1}, $\Theta^{(i)}$ with ${\bf X}^{(i)}$ satisfies a general enough weak dependence structure, that contains the independence as special case (${\bf X}^{(i)}$ has arbitrarily dependent components in all over the text). A frequent question, with both theoretical and practical value, is if the 
distribution of the random vector ${\bf X}$ belongs to some class of distributions, then does the distribution of the random vector $\Theta\,{\bf X}$ belong to the same class? 

We provide some sufficient conditions, for positive answer to this question for several distribution classes. Next, in subsection 4.2, we establish the asymptotic behavior of the scale mixture
\beam \label{eq.KP.1.2}
{\bf S}_n^{\Theta} = \sum_{i=1}^{n} \Theta^{(i)}\,{\bf X}^{(i)} \,,
\eeam
with $n \in \bbn$, as $\xto$, on a set $x\,A$, which can be interpreted as 'rare set', see in subsection 2.2 below. We obtain two results, where we establish the presence of the multivariate linear single big jump of the scale mixture sum in \eqref{eq.KP.1.2}. In both results, each pair $(\Theta^{(i)},\,{\bf X}^{(i)})$ contains this, already mentioned, dependence structure.

In Section 5, we provide the asymptotic behavior of the entrance probability of the discounted aggregate claims, into a rare-set $x\,A$, in the frame of a time-dependent multivariate risk model with common Poisson counting process and with constant interest rate. We give also the asymptotic expression of the finite-time ruin probability in several ruin-sets.

Finally, in Section 6, we present the lower bound for the vector type precise large deviations. In case of non-random sums, we assume weakly dependent vectors, with each of them to have arbitrarily dependent components, without any assumption for its distribution and we provide a uniform estimation for the lower bound. In case of random sums, under the additional assumption of weak equivalence of distributions of the random vectors, in the 'rare-set' of interest $x\,A$, and under the condition that the common counting process converges in probability to its mean value, as $t \to \infty$, we give uniform estimation for the lower bound. We also give an application of the last result to a multivariate  compound risk model driven by a quasi-renewal process.

\section{Heavy tails and preliminary results} \label{sec.KP.2}

In what follows, all the asymptotic relations hold as $\xto$, except otherwise stated, and the random vectors have support on the non-negative quadrant. 

For two positive one-dimensional functions $f,\,g$, we write $f(x) \sim c\,g(x)$, for some constant $c \in (0,\,\infty)$, if 
\beao
\lim \dfrac {f(x)}{g(x)} = c\,.
\eeao 
We denote by $f(x)=o[g(x)]$, if it holds 
\beao
\lim \dfrac {f(x)}{g(x)} = 0\,,
\eeao
and by  $f(x)=O[g(x)]$, if it holds 
\beao
\limsup \dfrac {f(x)}{g(x)} < \infty\,.
\eeao
Furthermore, we  denote by $f(x) \asymp g(x)$, if the relations $f(x)=O[g(x)]$ and $g(x)=O[f(x)]$ hold all together. Finally, we write $f(x) \lesssim g(x)$, (or, $f(x) \gtrsim g(x)$), if 
\beao
\limsup \dfrac{f(x)}{g(x)} \leq 1\,, \;\left({\text or},\;\liminf \dfrac{f(x)}{g(x)} \geq 1\,, \right)
\eeao
respectively.

Next, if the positive functions ${\bf f}$, ${\bf g}$, are $d$-variate, with $d\in \bbn$, then the corresponding notations hold for functions of the 
form ${\bf f}(x\,\bbb)$, ${\bf g}(x\,\bbb)$ for some set $\bbb \in \bbr^d$, for which it holds ${\bf 0} \not\in \bbb$. For example we write, 
${\bf f}(x\,\bbb)\sim c\,{\bf g}(x\,\bbb)$, for some $c \in (0,\,\infty)$, if it holds
\beao
\lim \dfrac {{\bf f}(x\,\bbb)}{{\bf g}(x\,\bbb)} = c\,.
\eeao
For any two $d$-dimensional vectors ${\bf x}$ and ${\bf y}$, all the operations are defined component-wise, as for example ${\bf x} \pm {\bf y}=(x_1\pm y_1,\,\ldots,\,x_d \pm y_d)^{\top}$. 
Further, the scalar product of $c\in (0,\,\infty)$ with ${\bf x}$ is defined as $c\,{\bf x}=(c\,x_1,\,\ldots,\,c\,x_d)^{\top}$. 

For two real numbers $x$, $y$, we denote by $x\vee y := \max\{x,\,y\}$ their maximum, by $x\wedge y := \min\{x,\,y\}$ their minimum, and by $\left\lfloor  x \right\rfloor$ we denote the integer part of $x$. By ${\bf 0}=(0,\,\ldots,\,0)^{\top}$ we denote the origin of the axes. By ${\bf 1}_{E}$ we denote the indicator function of the event $E$, and by 
$E^c$ the complementary set of $E$ with respect to the sample space $\Omega$. For any random variable $Z$, $Z \stackrel{d}{\sim} V$, indicates that random variable $Z$ has distribution $V$. For any distribution function $V$, we denote by $\overline{V}=1-V$ its tail and with $r_V$ its right endpoint 
of its support. For any two distributions $V_1$, $V_2$ we define by $V_1*V_2$ their convolution and by $V_1 V_2$ the distribution of the maximum, 
namely if $Z_1 \stackrel{d}{\sim} V_1$, $Z_2 \stackrel{d}{\sim} V_2$, with $Z_1$, $Z_2$ independent random variables, then $V_1*V_2(x) := \PP[Z_1+Z_2 \leq x]$ and 
$V_1 V_2(x):=\PP[Z_1 \vee Z_2 \leq x]$. With $V^{n*}$, with $n \in \bbn$, we denote the $n$-order convolution power of $V$. Finally we say that 
the distributions  $V_1$, $V_2$ have strongly equivalent tails if $\overline{V}_1(x) \sim c\,\overline{V}_2(x)$, for some $c \in (0,\,\infty)$ and 
they have weakly equivalent tails if $\overline{V}_1(x) \asymp \overline{V}_2(x)$.

\subsection{Heavy-tailed distributions} \label{sec.KP.2.1}

Here we refer some necessary preliminary definitions in relation with the classes of heavy-tailed distributions in one-dimensional case, as also 
in relation with some corresponding indexes, see for more details about properties and applications of these classes in \cite{embrechts:kluepellberg:mikosch:1997}, \cite{foss:korshunov:zachary:2013}, \cite{leipus:siaulys:konstantinides:2023}. Let us notice that for the sake of compactness of the paper, together with the multivariate distribution classes, we define the following one-dimensional distribution classes only for distributions $V$ with support on $\bbr_+=[0,\,\infty)$. In what follows 
in this subsection, we assume that $\bV(x)>0$, for all $x \in \bbr$.

We say that a distribution $V$ has heavy tail, symbolically $V \in \mathcal{K}$, if for any $\vep >0$, it holds
\beao
\int_{0-}^{\infty} e^{\vep\,x}\,V(dx)=\infty\,.
\eeao

Class $\mathcal{K}$ is a rather large one, and it contains many famous subclasses. A large subclass of $\mathcal{K}$ is the class  $\mathcal{L}$ of long tailed distributions, for which we say that $V \in \mathcal{L}$ if for any (or, equivalently, for some) $a>0$ it holds
\beao
\lim \dfrac{\overline{V}(x-a)}{\overline{V}(x)} = 1\,.
\eeao
Another important subclass is the class  $\mathcal{S}$ of subexponential distributions, for which we say that $V \in \mathcal{S}$ if for any (or, equivalently, for some) integer $n\geq 2$ it holds
\beao
\lim \dfrac{\overline{V^{n*}}(x)}{\overline{V}(x)} = n\,.
\eeao
The three classes $\mathcal{S},\,\mathcal{L},\,\mathcal{K}$ were introduced in \cite{chistyakov:1964}, where  the inclusions $\mathcal{S} \subsetneq \mathcal{L} \subsetneq \mathcal{K}$ were established.

Before proceeding to other subclasses of $\mathcal{K}$, we provide a concept related directly to their characterization. This is the Matuszewska indexes, 
see \cite{matuszewska:1964} for their introduction. We define the upper and lower Matuszewska indexes as follows. Let us denote
\beam \label{eq.KP.2.4}
\overline{V}_*(b):= \liminf \dfrac{\overline{V}(b\,x)}{\overline{V}(x)}\,, \qquad \overline{V}^*(b):= \limsup \dfrac{\overline{V}(b\,x)}{\overline{V}(x)}\,, 
\eeam 
for any $b>1$. Next, we define by 
\beam \label{eq.KP.2.3}
J_{V}^+:= -\lim_{b \to \infty} \dfrac{\log \overline{V}_*(b)}{\log b}\,, \qquad J_{V}^-:= -\lim_{b \to \infty} \dfrac{\log \overline{V}^*(b)}{\log b}\,, 
\eeam 
the upper and lower Matuszewska indexes. By relations \eqref{eq.KP.2.3} and  \eqref{eq.KP.2.4}, we see that for any distribution $V$, with $\bV(x)>0$, for any $x \in \bbr$, it holds $0\leq J_{V}^- \leq J_{V}^+ \leq \infty$. 

We say that a distribution $V$ belongs to the class of dominatedly varying distributions, symbolically $V \in \mathcal{D}$, if for any (or, equivalently, for some) $y \in (0,\,1)$, it holds
\beao
\limsup \dfrac{\overline{V}(y\,x)}{\overline{V}(x)}< \infty\,. 
\eeao

We say that a distribution $V$ belongs to the class of positively decreasing distributions, symbolically $V \in \mathcal{P_D}$, if for any (or, equivalently, for some) $b > 1$, it holds
\beao
\limsup \dfrac{\overline{V}(b\,x)}{\overline{V}(x)}< 1\,. 
\eeao 
These classes above were introduced by \cite{feller:1969}, and by \cite{dehaan:resnick:1984}, respectively, and they are symmetric in some sense. Indeed, we can see that the inclusion $V \in \mathcal{D}$ is equivalent to 
the inequality $J_V^+ < \infty$, while the inclusion $V \in \mathcal{P_D}$ is equivalent to the inequality $J_V^- >0$, see for example in \cite[Sec. 2.4]{leipus:siaulys:konstantinides:2023}. Furthermore, it is worth to mention that 
$\mathcal{D} \subsetneq \mathcal{K}$, while $\mathcal{P_D} \cap \mathcal{K} \neq \emptyset$ and $\mathcal{P_D} \setminus \mathcal{K} \neq \emptyset$, namely the class $\mathcal{P_D}$ contains both light-tailed and heavy-tailed distributions. Even more, we find in \cite{goldie:1978}, that $\mathcal{D} \not\subset \mathcal{S}$, $\mathcal{S} \not\subset \mathcal{D}$ and $\mathcal{D}\cap \mathcal{S} \equiv \mathcal{D}\cap \mathcal{L} \neq \emptyset$. Additionally, we obtain $\mathcal{D}\cap \mathcal{P_D} \subsetneq \mathcal{K}$. For more details on $\mathcal{P_D}$, see in \cite{bardoutsos:konstantinides:2011} and \cite{konstantinides:passalidis:2024d}.

In \cite{konstantinides:tang:tsitsiashvili:2002} we find the class of subexponential positively decreasing distributions, symbolically $\mathcal{A}:= \mathcal{S}\cap \mathcal{P_D}$, with motivation of the correction of a gap in 
\cite{kalashnikov:konstantinides:2000}, when it holds $J_V^-=0$. Furthermore, in several papers in risk theory was pointed out the importance to have strictly positive the lower Matuszewska index, as several times are needed assumptions for distributions from classes $\mathcal{A}$, $\mathcal{D}\cap \mathcal{A}$ and $\mathcal{D}\cap \mathcal{P_D}$, see for example in \cite{bardoutsos:konstantinides:2011}, \cite{dindiene:leipus:2015}, \cite{yang:zhang:jiang:cheng:2015} and \cite{chen:cheng:2024} among others.

Just recently, in \cite{konstantinides:passalidis:2024b} was introduced and studied the class of positively decreasing long tailed distributions, symbolically $\mathcal{T}:=\mathcal{L}\cap \mathcal{P_D}$. Later, in \cite{konstantinides:passalidis:2024d} the crucial role of this class for the closure property of class $\mathcal{A}$ with respect to convolution was shown.

Further, the class $\mathcal{C}$ of the consistently varying distributions is also included in the family of heavy-tailed distributions. We say that distribution $V \in \mathcal{C}$ if it holds
\beam \label{eq.KP.2.8}
\lim_{b\downarrow 1} \liminf \dfrac{\overline{V}(b\,x)}{\overline{V}(x)}= 1\,. 
\eeam 
It is well-known that $\mathcal{C} \subsetneq \mathcal{D}\cap \mathcal{L}$.

Finally, the class of regularly varying distributions is the most popular among the heavy-tailed distributions. We say that a distribution $V$ is regularly varying with index $\alpha \in [0,\,\infty)$, symbolically $V \in \mathcal{R}_{-\alpha}$, if it holds
\beao
\lim \dfrac{\overline{V}(t\,x)}{\overline{V}(x)}=t^{-\alpha}\,, 
\eeao 
for any $t>0$. It is well-known that if $V \in \mathcal{R}_{-\alpha}$ then $J_V^- = J_V^+=\alpha$ and $\bigcup_{0\leq \alpha <\infty} \mathcal{R}_{-\alpha} \subsetneq \mathcal{C}$. Furthermore, from the characterization of classes through Matuszewska indexes, it follows that
\beao
\bigcup_{0< \alpha <\infty} \mathcal{R}_{-\alpha} \subsetneq \mathcal{C} \cap \mathcal{P_D}\,. 
\eeao
For further discussions and treatments about regular variation, we refer to 
\cite{bingham:goldie:teugels:1987}. It is worth reminding the class of rapidly 
varying distributions $\mathcal{R}_{-\infty}$, that characterize distributions $V$ that satisfy relation
\beao
\lim \dfrac{\overline{V}(b\,x)}{\overline{V}(x)}=0
\eeao
for any $b>1$, obviously satisfy the inclusion $\mathcal{R}_{-\infty} \subsetneq \mathcal{P_D}$.

\subsection{Heavy-tailed random vectors} \label{sec.KP.2.2}

Let us mention some well-known definitions of multivariate distributions with heavy tails and next we introduce the class of multivariate positively decreasing distributions, as well as its intersection with some other classes. We note that the most popular multivariate heavy tailed class is the multivariate regular variation $MRV$, introduced in \cite{dehaan:resnick:1981}. If ${\bf X} \stackrel{d}{\sim} F$ we say that $F$ belongs to the (standard) $MRV$, if there exists a Radon measure $\mu$, that is non-degenerate at zero, and an one-dimensional distribution 
$V \in \mathcal{R}_{-\alpha}$, with $\alpha \in (0,\,\infty)$, such that it holds
\beam \label{eq.KP.2.11}
\lim \dfrac 1{\overline{V}(x)} \PP[{\bf X} \in x\,\bbb] = \mu(\bbb)\,, 
\eeam 
for any Borel, $\mu$-continuous set $\bbb \in [0,\,\infty]^d $, that is bounded 
away from ${\bf 0}$. We denote this multivariate regular variation distribution as $F\in MRV(\alpha,\,V,\,\mu)$. Even more, the Radon measure in \eqref{eq.KP.2.11} has the property of positive homogeneity (of order $-\alpha$), in the sense that for any Borel set $\bbb \in  [0,\,\infty]^d$ (bounded away from ${\bf 0}$) and any constant $c>0$, it holds
\beam \label{eq.KP.2.12}
\mu(c^{1/\alpha}\,\bbb)=\dfrac 1{c} \mu(\bbb) \,, 
\eeam 
see for further material in \cite{resnick:2007}, \cite{mikosch:wintenberger:2024}. However, although the standard $MRV$ is well studied, it is not the case of the non-standard $MRV$. 

Now we need the set family $\mathscr{R}:=\left\{A\subsetneq \bbr^d\;:\;A \;\text{open},\;\text{increasing},\;A^c\;\text{convex},\; {\bf 0} \notin \overline{A} \right\}$, where by $\overline{A}$ we denote the closure of set $A$. A set $A$ is called increasing, if for any ${\bf x} \in A$ and any ${\bf y}\in \bbr_+^d$, we obtain that ${\bf x}+{\bf y} \in A$. From \cite[Lem. 4.3 (c)]{samorodnitsky:sun:2016} we find that for any $A \in \mathscr{R}$, there exists a set $I_A \in \bbr^{d}$, such that it holds
\beam \label{eq.KP.2.14}
A=\left\{ {\bf x} \in \bbr^{d}\;:\;{\bf p}^{\top}\,{\bf x}=\sum_{i=1}^d\,p_i\,x_i >1\,, \;\exists \;{\bf p} \in I_A  \right\} \,.
\eeam 
Let us notice that if $A \in \mathscr{R}$, then for any constant $c>0$, we obtain $c\,A \in \mathscr{R}$, namely the set family $\mathscr{R}$ represent a cone with respect to multiplication with positive scalars. Through this property and relation \eqref{eq.KP.2.14}, in \cite[Lem. 4.5]{samorodnitsky:sun:2016} was proved that for the random vector ${\bf X} \stackrel{d}{\sim} F$, the random variable
\beam \label{eq.KP.2.15}
Y_A:=\sup \left\{u\;:\; {\bf X} \in u\,A \right\} \,
\eeam   
follows a proper distribution $F_A$, whose tail is defined by
\beam \label{eq.KP.2.16}
\overline{F}_A(x)=\PP\left[{\bf X} \in x\,A \right] = \PP\left[\sup_{{\bf p} \in I_A} {\bf p}^{\top}\,{\bf X}> x \right] \,, \quad x>0\,.
\eeam  
With the help of relation \eqref{eq.KP.2.16}, in \cite{samorodnitsky:sun:2016}  the multivariate subexponentiality was defined. Let assume that $A\in \mathscr{R}$ is a fixed set. We say that the vector ${\bf X} \stackrel{d}{\sim} F$ has multivariate subexponential distribution on $A$, symbolically $F \in \mathcal{S}_A$, if $F_A \in \mathcal{S}$.

Next, from \cite{konstantinides:passalidis:2024g}, we obtain the definition of multivariate dominatedly varying distributions on $A$, symbolically $F\in \mathcal{D}_A$, if $F_A \in \mathcal{D}$, the definition of multivariate long tailed  distributions on $A$, symbolically $F \in \mathcal{L}_A$, if $F_A \in \mathcal{L}$, the definition of consistently varying distributions on $A$, symbolically $F\in \mathcal{C}_A$, if $F_A \in \mathcal{C}$ and the definition of multivariate dominatedly varying, long tailed distributions on $A$, symbolically $F\in (\mathcal{D}\cap\mathcal{L})_A$, if $F_A \in \mathcal{D}\cap \mathcal{L}$. For all these distribution classes we denote 
\beam \label{eq.KP.2.17}
\mathcal{B}_{\mathscr{R}}:=\bigcap_{A\in \mathscr{R} } \,\mathcal{B}_A \,,
\eeam 
for $\mathcal{B} \in \{\mathcal{C},\,\mathcal{D},\,\mathcal{D}\cap \mathcal{L},\,\mathcal{S},\,\mathcal{L} \}$, as their multivariate analogue on whole set family $\mathscr{R}$. Now, from relations \eqref{eq.KP.2.14}, \eqref{eq.KP.2.16} 
(see also, \cite[Rem. 4.1]{samorodnitsky:sun:2016}) and the definition of previous multivariate distribution classes, according to description in \cite{samorodnitsky:sun:2016} and \cite{konstantinides:passalidis:2024g}, if ${\bf X} \stackrel{d}{\sim} F \in \mathcal{B}_{\mathscr{R}}$ then the distribution of all the non-degenerate at zero and non-negative  linear combinations, namely the distributions of $\sum_{i=1}^d l_i\,X_i$, belong to class $\mathcal{B}$. Several well-established distribution classes, as the class of multivariate stable and multivariate infinite divisible distributions dispose this property. Let us recall that if $F \in \mathcal{B}_A$ for some  $A \in \mathscr{R}$ but $F \not\in \mathcal{B}_\mathscr{R}$, then the closure property of linear combinations is not necessarily valid.

Now, we define the class of multivariate positively decreasing distributions, as also its intersection with other multivariate classes.
\bde \label{def.KP.2.1}
Let some fixed set $A \in \mathscr{R}$ and some random vector ${\bf X} \stackrel{d}{\sim} F$. We say that $F$ belongs to the class of multivariate positively decreasing distributions on $A$, symbolically $F \in \mathcal{P_D}_A$, if it holds that $F_A \in \mathcal{P_D}$.

We say that $F$ belongs to the class of multivariate long-tailed positively decreasing distributions on $A$, symbolically $F \in \mathcal{T}_A$, if it holds that $F_A \in \mathcal{T}$.

We say that $F$ belongs to the class of multivariate subexponential positively decreasing distributions on $A$, symbolically $F \in \mathcal{A}_A$, if it holds that $F_A \in \mathcal{A}$.

We say that $F$ belongs to the class of multivariate dominatedly varying, positively decreasing distributions on $A$, symbolically $F \in (\mathcal{D}\cap \mathcal{P_D})_A$, if it holds that $F_A \in \mathcal{D}\cap \mathcal{P_D}$.

We say that $F$ belongs to the class of multivariate dominatedly varying, subexponential, positively decreasing distributions on $A$, symbolically $F \in (\mathcal{D}\cap \mathcal{A})_A$, if it holds that $F_A \in \mathcal{D}\cap \mathcal{A}$.

We say that $F$ belongs to the class of multivariate consistently varying, positively decreasing distributions on $A$, symbolically $F \in (\mathcal{C}\cap \mathcal{P_D})_A$, if it holds that $F_A \in \mathcal{C}\cap \mathcal{P_D}$.
\ede 

\bre \label{rem.KP.2.1}
In the previous definitions, we can easily see that $(\mathcal{D}\cap \mathcal{T})_A:= (\mathcal{D}\cap \mathcal{A})_A$, in the sense that $F\in (\mathcal{D}\cap \mathcal{T})_A$ if it holds $F_A \in \mathcal{D}\cap \mathcal{T}$, because of the equality $\mathcal{D}\cap \mathcal{S} \equiv \mathcal{D}\cap \mathcal{L}$. Furthermore, we keep the notation of \eqref{eq.KP.2.17} for all the classes $\mathcal{B} \in \{\mathcal{P_D},\,\mathcal{T},\,\mathcal{A},\, \mathcal{D}\cap \mathcal{P_D},\,\mathcal{D}\cap \mathcal{A},\,\mathcal{C} \cap \mathcal{P_D}\}$ and for anyone of them, directly by Definition \ref{def.KP.2.1}, relation \eqref{eq.KP.2.16}(and \cite[Rem. 4.1]{samorodnitsky:sun:2016}), we obtain that if $F \in \mathcal{B}_\mathscr{R}$, then the distribution of 
\beao
\sum_{i=1}^d l_i\,X_i\,, 
\eeao
for non-negative and non-degenerate at zero, linear combinations, belongs to $\mathcal{B}$. Again, this is not necessarily valid if $F  \in \mathcal{B}_A$ (for some $A \in \mathscr{R}$) but $F \not\in \mathcal{B}_\mathscr{R}$. For some examples for classes $\mathcal{B}_\mathscr{R},\,\mathcal{B}_A$, with $\mathcal{B} \in \{\mathcal{C} \cap \mathcal{P_D},\,\mathcal{D}\cap \mathcal{A},\,\mathcal{A},\,\mathcal{S}\}$, that escape from the $MRV$ frame, 
we refer to \cite[Sec. 4]{samorodnitsky:sun:2016} and \cite[Sec. 4]{konstantinides:liu:passalidis:2026}.
\ere

\section{Convolution closure and related properties} \label{sec.KP.3}

Here, we provide some properties of the new classes introduced in Definition \ref{def.KP.2.1}. The next result gives an inclusion relation.

\bpr \label{pr.KP.2.1}
Let $A \in \mathscr{R}$, and ${\bf X} \stackrel{d}{\sim} F \in MRV(\alpha,\,V,\,\mu)$, with $\alpha \in (0,\,\infty)$. If $\mu(A) >0$, then $F\in (\mathcal{C}\cap \mathcal{P_D})_{A}$.
\epr 

\pr~
In the proof of \cite[Prop. 4.14]{samorodnitsky:sun:2016} was established that for any $A \in \mathscr{R}$ it holds $\mu(\partial A)=0$ and $\mu(A)< \infty$, (because, ${\bf 0} \not\in \overline{A}$). So by \eqref{eq.KP.2.11}, \eqref{eq.KP.2.16} and the assumption $\mu(A)>0$, we obtain
\beao
\PP\left[Y_A >x\right] \sim \mu(A)\,\overline{V}(x)\,.
\eeao
So, since $V \in \mathcal{R}_{-\alpha}$, with $\alpha \in (0,\,\infty)$ from the closure property of the regular variation with respect to strong equivalence, see \cite[Prop. 3.3(i)]{leipus:siaulys:konstantinides:2023}, follows that $F_A \in \mathcal{R}_{-\alpha}$ and since $\mathcal{R}_{-\alpha} \subsetneq \mathcal{C} \cap \mathcal{P_D}$, we find out that $F_A \in \mathcal{C} \cap \mathcal{P_D}$. Therefore we obtain $F\in (\mathcal{C}\cap \mathcal{P_D})_A$. 
~\halmos

\bre \label{rem.KP.2.2}
From the previous proposition (in combination with the definitions of our classes) it seems that the one-dimensional inclusions with respect to the order of the classes remains the same in the multivariate set up. Hence, for any $A \in \mathscr{R}$, such that $\mu(A) >0$ it holds 
\beam \label{eq.KP.2.19}
\bigcup_{0<\alpha < \infty} MRV(\alpha,\,V,\,\mu) \subsetneq \mathcal{C}_{A}  \subsetneq (\mathcal{D}\cap \mathcal{L})_{A} \subsetneq \mathcal{S}_{A} \subsetneq \mathcal{L}_{A}\,,
\eeam
and
\beao
\bigcup_{0<\alpha < \infty} MRV(\alpha,\,V,\,\mu) \subsetneq (\mathcal{C} \cap \mathcal{P_D})_A \subsetneq  (\mathcal{D}\cap \mathcal{A})_{A} \subsetneq (\mathcal{D} \cap \mathcal{P_D})_{A}\,.
\eeao
The previous inclusions are valid for the corresponding classes $\mathcal{B}_{\mathscr{R}}$, when $\mu(A)>0$, for any $A \in \mathscr{R}$. 
\ere

\bre \label{rem.KP.3.A}
In fact, from the proof of last proposition we find that if $F \in MRV(\alpha,\,V,\,\mu)$, then for any $A \in \mathscr{R}$, such that $\mu(A)>0$, follows that $F_A \in \mathcal{R}_{-\alpha}$. The inverse is not true. In fact, if $\mu(A) =0$, then the $MRV$ and $F_A \in \mathcal{R}_{-\alpha}$, are incomparable as is implied by \cite[Exam. 2.1]{chen:passalidis:2026}. Determining a concrete $A \in \mathscr{R}$, the classes of \eqref{eq.KP.2.19}, can provide a more general and more useful tool for the actuarial mathematics. These classes can contain cases of marginal distributions, that have not necessarily weak-equivalent tails, and can describe precisely the decay rate of the ruin probability or related quantities (in such cases, the using of $MRV$ gives frequently zero Radon measure). We refer to \cite[Exam. 4.1 - 4.4]{konstantinides:liu:passalidis:2026}, for some cases, that permit inhomogeneity of the marginals.   
\ere

A key closure property in one-dimensional heavy-tailed distributions is the closure with respect to strong tail equivalence, that means if $V_1 \in \mathcal{B}$ and  $\overline{V}_2(x) \sim c\,\overline{V}_1(x)$ for some $c\in (0,\,\infty)$, then it holds $V_2 \in \mathcal{B}$. This property is satisfied for class $\mathcal{P_D}$, since if $V_1 \in \mathcal{P_D}$ and $\overline{V}_2(x) \sim c\,\overline{V}_1(x)$ for some $c\in (0,\,\infty)$, then for any $v>1$ it holds
\beam \label{eq.KP.3.A}
\limsup \dfrac{\overline{V}_2(v\,x)}{\overline{V}_2(x)} = \limsup \dfrac{c\,\overline{V}_1(v\,x)}{c\,\overline{V}_1(x)} < 1\,.
\eeam

In the next proposition we present a kind of closure property with respect to strong tail equivalence for the multivariate classes of Definition \ref{def.KP.2.1}. In what follows we consider that $F_i(x\,A) = \PP\left[{\bf X}^{(i)} \in x\,A\right]$, for any $i=1,\,\ldots,\,n$ with $n\in \bbn$, and $Y_A^{(i)}:=\sup \left\{u\;:\;{\bf X}^{(i)} \in u\,A\right\}$, for any $i=1,\,\ldots,\,n$, with 
$Y_A^{(i)} \stackrel{d}{\sim} F_A^{(i)}$.

\bpr \label{pr.KP.2.2}
Let $A \in \mathscr{R}$ be a fixed set and let random vectors ${\bf X}^{(1)}$, ${\bf X}^{(2)}$ be arbitrary dependent, with distributions $F_1$,\, $F_2$, respectively. If $F_1 \in \mathcal{B}_A$ and 
\beam \label{eq.KP.2.21}
\lim \dfrac{F_2(x\,A)}{F_1(x\,A)}=c\,,
\eeam
for some constant $c \in (0,\,\infty)$, then $F_2 \in \mathcal{B}_A$, with $\mathcal{B} \in \{\mathcal{P_D},\,\mathcal{T},\,\mathcal{A},\, \mathcal{D}\cap \mathcal{P_D},\,\mathcal{D}\cap \mathcal{A},\,\mathcal{C} \cap \mathcal{P_D}\}$.
\epr

\pr~
At first, for inclusion $F_2 \in \mathcal{B}_A$, it is sufficient to show $F_A^{(2)} \in \mathcal{B}$. From relation \eqref{eq.KP.2.21}, we find directly that 
\beao
\overline{F_A^{(2)}}(x) \sim c\, \overline{F_A^{(1)}}(x)\,. 
\eeao
When $\mathcal{B}=\mathcal{P_D}$, by the closure property of class $\mathcal{P_D}$ with respect to strong equivalence, see relation \eqref{eq.KP.3.A}, we obtain that $F_A^{(2)} \in \mathcal{P_D}$, since $F_A^{(1)} \in \mathcal{P_D}$, which implies $F_2 \in \mathcal{P_D}_A$. For the rest subclasses, the result follows from the closure property of $\mathcal{P_D}_A$ and \cite[Prop. 4.12(a)]{samorodnitsky:sun:2016} and \cite[Prop. 2.2]{konstantinides:passalidis:2024g}.
~\halmos

When one of the two multivariate distributions belongs to class $\mathcal{S}_A$ (or, to class  $\mathcal{A}_A$), and the other to class $\mathcal{L}_A$, the following proposition provides a closure property of $\mathcal{S}_A$, or of $\mathcal{A}_A$, with respect to weak tail-equivalence. The corresponding one-dimensional result for class $\mathcal{S}$ can be found in \cite{kluppelberg:1988}. It is also easy to verify that if $\bV_1(x) \asymp \bV_2(x)$, then we obtain $J_{V_1}^{-}=J_{V_2}^{-}$. Hence, class $\mathcal{P_D}$ is closed with respect to weak tail-equivalence.

\bpr \label{prop.KP.3.A*}
Let $A \in \mathscr{R}$ be some fixed set and ${\bf X}^{(1)}$, ${\bf X}^{(2)}$ be arbitrarily dependent random vectors with distributions $F_1,\,F_2$, respectively. If $F_1  \in \mathcal{S}_A$ (or, in $\mathcal{A}_A$) and $F_2  \in \mathcal{L}_A$, and additionally the relationship
\beam \label{eq.KP.3.a}
F_1(x\,A) \asymp F_2(x\,A)\,,
\eeam 
is valid, then it holds $F_2 \in \mathcal{S}_A$ (respectively, in $\mathcal{A}_A$).
\epr

\pr~
Let start with class $\mathcal{S}_A$. From relation \eqref{eq.KP.3.a}, is implied
\beam \label{eq.KP.3.b}
\bF_A^{(1)}(x) \asymp \bF_A^{(2)}(x)\,.
\eeam 
From the fact that $F_A^{(1)} \in \mathcal{S}$ and $F_A^{(2)} \in \mathcal{L}$, and relation \eqref{eq.KP.3.b} is valid, through the \cite[Th. 2.1(a)]{kluppelberg:1988} or, the \cite[Prop. 3.13(ii)]{leipus:siaulys:konstantinides:2023}, we obtain that $F_A^{(2)} \in \mathcal{S}$, hence we get $F_2 \in \mathcal{S}_A$.

Now we study the other  part, with class $\mathcal{A}_A$. Since, we know that $\mathcal{A} \subsetneq \mathcal{S}$, and we assumed \eqref{eq.KP.3.a}, from previous arguments we find that $F_A^{(2)} \in \mathcal{S}$. Because of $F_A^{(1)} \in \mathcal{P_D}$, from \eqref{eq.KP.3.b} we have $F_A^{(2)} \in \mathcal{P_D}$. So, we get $F_A^{(2)} \in \mathcal{A}$, from where we find $F_2 \in \mathcal{A}_A$.
~\halmos

For two independent random vectors we define their  convolution by
\beao
F_1*F_2(x\,A)=\PP\left[{\bf X}^{(1)}+{\bf X}^{(2)} \in x\,A\right]\,.
\eeao 
So for some class $\mathcal{B}_A$ we say that it is closed with respect to convolution if $F_1*F_2 \in \mathcal{B}_A$, namely if the random variable 
\beao
Y_A^*:=\sup\left\{u\;:\;{\bf X}^{(1)}+{\bf X}^{(2)} \in u\,A\right\} \stackrel{d}{\sim} F_A^*\,,
\eeao 
satisfies the inclusion $F_A^* \in \mathcal{B}$. Let us notice as shown in \cite[Prop. 2.4]{konstantinides:passalidis:2024g} for arbitrarily dependent ${\bf X}^{(1)},\,\ldots,\,{\bf X}^{(n)}$, it holds for any $A \in \mathscr{R}$ that 
\beam \label{eq.KP.2.22}
\PP[{\bf X}^{(1)}+\,\cdots\,+{\bf X}^{(n)} \in x\,A] \leq \PP\left[ \sum_{i=1}^n Y_A^{(i)} >x\right]\,, \quad x\geq 0\,.
\eeam
The following result gives a closure property with respect to convolution for subclass $(\mathcal{D}\cap \mathcal{A})_A$. We observe that in one-dimensional case for $A=(1,\,\infty)$, the following relation coincides with \cite[Prop. 3.2]{bardoutsos:konstantinides:2011}.

\bpr \label{pr.KP.2.3}
Let $A \in \mathscr{R}$ be some fixed set and  ${ \bf X}^{(1)},\,{\bf X}^{(2)}$ be independent random vectors with distributions $F_1,\,F_2 \in \mathcal{B}_A=\{(\mathcal{C} \cap \mathcal{P_D})_A,\, (\mathcal{D}\cap \mathcal{A})_A \}$ respectively. Then it holds
\beam \label{eq.KP.2.23}
F_1*F_2(x\,A) \sim F_1(x\,A)+F_2(x\,A)\,,
\eeam
and further $F_1*F_2  \in \mathcal{B}_A$.
\epr

\pr~
From \cite[Th. 3.3(ii)]{konstantinides:passalidis:2024g} and the fact that $\mathcal{B}_A \subsetneq (\mathcal{D}\cap \mathcal{L})_A$, we obtain that relation \eqref{eq.KP.2.23} is true. Further, by \cite[Th. 3.3(ii), (iii)]{konstantinides:passalidis:2024g}, we obtain that $F_1*F_2  \in (\mathcal{D}\cap \mathcal{L})_A$ (or $F_1*F_2  \in \mathcal{C}_A$, respectively). For any $v>1$ it holds
\beao
\limsup \dfrac{\PP\left[Y_A^* > v\,x\right]}{\PP\left[Y_A^* > x\right]}&=& \limsup \dfrac{F_1*F_2(v\,x\,A)}{F_1*F_2(x\,A)}= \limsup \dfrac{F_1(v\,x\,A)+F_2(v\,x\,A)}{F_1(x\,A)+F_2(x\,A)}\\[2mm]
&=& \limsup \dfrac{\PP\left[ Y_A^{(1)} >v\,x\right]+\PP\left[  Y_A^{(2)} >v\,x\right]}{\PP\left[ Y_A^{(1)} >x\right]+\PP\left[  Y_A^{(2)} >x\right]} \\[2mm]
&\leq& \max_{1\leq i \leq 2}\left\{ \limsup \dfrac{\PP\left[ Y_A^{(i)} >v\,x\right]}{\PP\left[ Y_A^{(i)} >x\right]}\right\}< 1\,,
\eeao
where in the second step we used relation \eqref{eq.KP.2.23} and in the last step we used that $F_A^{(i)} \in \mathcal{P_D}$, for $i=1,\,2$. Thus, we find $F_1*F_2  \in \mathcal{P_D}_A$ and in combination with previous we obtain that $F_1*F_2  \in \mathcal{B}_A$.
~\halmos

\bre \label{rem.KP.2.3}
One of the famous problems in closure properties, is the closure of the subexponential class with respect to convolution. The well-known counter-example in \cite{leslie:1989} shows that the subexponential class is not closed with respect to convolution. A partial answer for this kind of closure, was presented in \cite{embrechts:goldie:1980}, while in \cite{leipus:siaulys:2020} was extended through equivalent conditions to the closure with respect to convolution, given that the distributions  belong to class $\mathcal{L}$. Inspired by that paper \cite[Th. 3.2]{konstantinides:passalidis:2024d} has shown similar results for class $\mathcal{A}$. More precisely, if $Y_A^{(1)}$,  $Y_A^{(2)}$ are independent random variables with distributions $F_A^{(1)},\,F_A^{(2)} \in \mathcal{T}$, it is established that the following statements are equivalent
\begin{enumerate}
\item
$F_A^{(1)}*F_A^{(2)}  \in \mathcal{A}$,
\item
$F_A^{(1)}F_A^{(2)}  \in \mathcal{A}$,
\item
$p\,F_A^{(1)} + (1-p)\,F_A^{(2)}  \in \mathcal{A}$, for any (or, equivalently, for some) $p \in (0,\,1)$.
\end{enumerate}
Furthermore, any of these equivalent statements implies that
\beam \label{eq.KP.2.24}
\overline{F_A^{(1)}*F_A^{(2)}}(x) \sim \overline{F_A^{(1)}}(x)+\overline{F_A^{(2)}}(x)\,.
\eeam
\ere
It  was shown in \cite[Cor. 3.1]{konstantinides:passalidis:2024d} that if we restrict the initial condition for the distribution classes of $F_A^{(1)}$, $F_A^{(2)}$ to be instead from $\mathcal{T}$, only just from $\mathcal{A}$, then the points  $(1)$ - $(3)$ and relation \eqref{eq.KP.2.24} are equivalent. Now, we can present necessary and sufficient conditions for the closure property of class $\mathcal{A}_A$.

\bth \label{lem.KP.2.1}
Let $A \in \mathscr{R}$ be some fixed set. We assume that ${\bf X}^{(1)}$, ${\bf X}^{(2)}$ are independent random vectors with distributions $F_1,\,F_2 \in \mathcal{T}_A$, respectively. Then it holds $F_A^{(1)}*F_A^{(2)}  \in \mathcal{A}$, if and only if $F_1*F_2  \in \mathcal{A}_A$. 
\ethe

\pr~
$(\Rightarrow)$. Assume that $F_A^{(1)}*F_A^{(2)}  \in \mathcal{A}$. At first, because of  the inclusion $ \mathcal{T}_A \subsetneq  \mathcal{L}_A$, from \cite[Th. 3.4]{konstantinides:passalidis:2024g}, we obtain that $F_1*F_2  \in \mathcal{S}_A$. Thus, it remains to show that $F_1*F_2  \in \mathcal{P_D}_A$. Because of the fact that ${\bf X}^{(1)}$, ${\bf X}^{(2)}$ are non-negative independent random vectors, and the $x\,A$ is increasing set, through Bonferroni inequality we get
\beam \label{eq.KP.2.25} \notag
&&\PP\left[{\bf X}^{(1)}+ {\bf X}^{(2)} \in x\,A\right] \geq \PP\left[\bigcup_{i=1}^2 \{{\bf X}^{(i)} \in x\,A\}\right] \\[2mm]
&&\geq \sum_{i=1}^2 \PP\left[{\bf X}^{(i)} \in x\,A\right] -  \sum_{1\leq i< j \leq 2} \PP\left[{\bf X}^{(i)} \in x\,A,\; {\bf X}^{(j)} \in x\,A\right]  \\[2mm] \notag
&&=\PP\left[ Y_A^{(1)} >x\right]+\PP\left[  Y_A^{(2)} >x\right] -\PP\left[ Y_A^{(1)} >x\right]\,\PP\left[  Y_A^{(2)} >x\right]  \\[2mm] \notag
&&\sim \PP\left[ Y_A^{(1)} >x\right]+\PP\left[  Y_A^{(2)} >x\right] -o\left(\PP\left[ Y_A^{(1)} >x\right]\right)\sim F_1(x\,A) + F_2(x\,A)\,.
\eeam
Hence, by relations \eqref{eq.KP.2.22} and \eqref{eq.KP.2.25}, for any $v>1$ we find
\beam \label{eq.KP.2.26} \notag
&&\limsup \dfrac{\PP\left[Y_A^{*} > v\,x\right]}{\PP\left[Y_A^{*} > x\right]} = \limsup \dfrac{\PP\left[{\bf X}^{(1)} +{\bf X}^{(2)} \in v\,x\,A\right]}{\PP\left[{\bf X}^{(1)} +{\bf X}^{(2)} \in x\,A\right]} \leq \limsup \dfrac{\PP\left[Y_A^{(1)}+ Y_A^{(2)}> v\,x\right]}{F_1(x\,A) + F_2(x\,A)}\\[2mm] 
&&= \limsup \dfrac{\overline{F_A^{(1)}}(v\,x) +\overline{F_A^{(2)}}(v\,x) }{\overline{F_A^{(1)}}(x) +\overline{F_A^{(2)}}(x) } \leq \max_{1\leq i \leq 2} \left\{ \limsup \dfrac{\overline{F_A^{(i)}}(v\,x)}{\overline{F_A^{(i)}}(x)}\right\}<1\,,
\eeam
where in the third step, we used relation \eqref{eq.KP.2.24}, which is valid because of the inclusions $F_A^{(1)},\,F_A^{(2)} \in \mathcal{T}$ and $F_A^{(1)}*F_A^{(2)} \in \mathcal{A}$, see Remark \ref{rem.KP.2.3}. Therefore, by \eqref{eq.KP.2.26} we find $F_1*F_2  \in \mathcal{P_D}_A$ and in combination with the previous statements we conclude that $F_1*F_2  \in \mathcal{A}_A$. 

$(\Leftarrow)$. Let $F_1*F_2  \in \mathcal{A}_A$. From \cite[Th. 3.4]{konstantinides:passalidis:2024g}, taking into account that 
$F_1,\,F_2 \in \mathcal{T}_A \subsetneq \mathcal{L}_A\,,\;\mathcal{A}_A \subsetneq \mathcal{S}_A$, we obtain  $F_A^{(1)}*F_A^{(2)} \in \mathcal{S}$. 
Hence, by $F_A^{(1)},\,F_A^{(2)} \in \mathcal{T} \subsetneq \mathcal{L}$ through \cite[Th. 1.1]{leipus:siaulys:2020} we obtain the validity of \eqref{eq.KP.2.24}. Therefore, 
for any $v>1$, it holds
\beao
\limsup \dfrac{\overline{F_A^{(1)}*F_A^{(2)}}(v\,x) }{\overline{F_A^{(1)}*F_A^{(2)}}(x)} &=& \limsup \dfrac{\bF_A^{(1)}(v\,x) + \bF_A^{(2)}(v\,x)}{\bF_A^{(1)}(x) +\bF_A^{(2)}(x)} \\[2mm]
&\leq&  \max_{1\leq i \leq 2} \left\{ \limsup \dfrac{\bF_A^{(i)}(v\,x)}{\bF_A^{(i)}(x)}\right\}<1\,,
\eeao
that implies $F_A^{(1)}*F_A^{(2)} \in \mathcal{P_D}$, and in combination with previous arguments, we conclude that $F_A^{(1)}*F_A^{(2)} \in \mathcal{A}$.
~\halmos

Now we study the closure property of classes $\mathcal{S}_A$, $\mathcal{A}_A$ with respect to convolution roots. Let us note that for some fixed set 
$A \in \mathscr{R}$, we say  that the class $\mathcal{B}$ is closed with respect to convolution roots, if from the $F^{n*} \in \mathcal{B}_A$, for some 
$n\geq 2$, then it follows that  $F\in \mathcal{B}_A$.

For one-dimensional heavy tailed distribution classes, there are many papers devoted to closure and non-closure properties, with respect to the 
convolution roots, see for example \cite{embrechts:goldie:veraverbeke:1979}, \cite{shimura:watanabe:2005}, \cite{watanabe:2008}, \cite{xu:foss:wang:2015}, 
\cite{cui:wang:xu:2022}, among others. Here, we provide a conditional positive answer for the closure properties of  classes 
$\mathcal{S}_A$, $\mathcal{A}_A$ with respect to convolution roots. Before answering to this issue, we need a preliminary lemma.

\ble \label{lem.KP.3.a}
Let $A \in \mathscr{R}$ be some fixed set. If either $F^{n*}_A \in \mathcal{S}$ or $F^{n*}\in \mathcal{S}_A$, for some integer $n\geq 2$, then it 
holds $\overline{F_A^{n*}}(x) \asymp F^{n*}(x\,A)$.
\ele

\pr~
Let $n\geq 2$ be some integer. At first we consider $F^{n*}_A \in \mathcal{S}$. Then for the random variable
$Y'_A :=\sup \{u\;:\;{\bf X}^{(1)} +\cdots + {\bf X}^{(n)} \in u\,A \} \stackrel{d}{\sim} F_A'$, we obtain from \cite[Lem. 4.9]{samorodnitsky:sun:2016} the relation
\beam \label{eq.KP.3.27}
\overline{F_A'}(x) = \PP[{\bf X}^{(1)} +\cdots + {\bf X}^{(n)} \in x\,A] \leq \PP[Y_A^{(1)} +\cdots + Y_A^{(n)}> x]=\overline{F_A^{n*}}(x)\,,
\eeam
for $Y_A^{(1)},\,\ldots,\, Y_A^{(n)}$ independent and identically distributed random variables, with common distribution $F_A$. From the other side, since the class $ \mathcal{S}$ is closed with respect to convolution roots, see \cite[Th. 2]{embrechts:goldie:veraverbeke:1979}, we find  $F_A \in \mathcal{S}$. Hence,
\beam \label{eq.KP.3.28}
\limsup \dfrac{\overline{F_A^{n*}}(x)}{F^{n*}(x\,A)} = \limsup \dfrac{n\,\bF_A(x)}{\overline{F_A'}(x) } \leq n < \infty\,,
\eeam
where in the pre-last step we used the non-negativeness of ${\bf X}^{(1)},\,\ldots,
\,{\bf X}^{(n)}$, remember relations \eqref{eq.KP.2.15} and \eqref{eq.KP.2.16}. So by relations \eqref{eq.KP.3.27} and \eqref{eq.KP.3.28}, we have the asymptotic relation $ \overline{F_A^{n*}}(x)\asymp F^{n*}(x\,A)$.

At second we assume $F^{n*}\in \mathcal{S}_A$. Now, the random variable  $Y'_A$, defined above, follows the distribution $F_A' \in \mathcal{S}$. Hence, if the  $Y_A'^{(1)},\,\ldots,\, Y_A'^{(n)}$ are independent and identically distributed random variables with distribution $F_A'$, from the non-negativeness of ${\bf X}^{(1)},\,\ldots,\,{\bf X}^{(n)}$ we obtain
\beam \label{eq.KP.3.29} \notag
&&\limsup \dfrac{\overline{F_A^{n*}}(x)}{F^{n*}(x\,A)} = \limsup \dfrac{\PP[Y_A^{(1)}+\cdots+ Y_A^{(n)}>x]}{\overline{F_A'}(x) } \\[2mm]
&&\leq \limsup \dfrac{\PP[Y_A'^{(1)}+\cdots+ Y_A'^{(n)}>x]}{\overline{F_A'}(x) }= n \,,
\eeam
By combination of relations \eqref{eq.KP.3.27} and \eqref{eq.KP.3.29} we find the desired result.
~\halmos

\bre \label{rem.KP.3.3b}
We can observe that either $F^{n*}_A \in \mathcal{S}$ or $F^{n*}\in \mathcal{S}_A$ imply the inequalities
\beao
1\leq \liminf \dfrac{\overline{F_A^{n*}}(x)}{F^{n*}(x\,A)} \leq \limsup \dfrac{\overline{F_A^{n*}}(x)}{F^{n*}(x\,A)} \leq n \,.
\eeao
\ere

\bth \label{th.KP.3.1b}
Let $A \in \mathscr{R}$ be some fixed set. 
\begin{enumerate}
\item[(i)]
If $F^{n*}\in \mathcal{S}_A$, for some integer $n\geq 2$, and $F \in \mathcal{L}_A$, then  $F \in \mathcal{S}_A $.
\item[(ii)]
If $F^{n*}\in \mathcal{A}_A$, for some integer $n\geq 2$, and $F \in \mathcal{L}_A$, then  $F \in \mathcal{A}_A $.
\end{enumerate}
\ethe

\pr~
\begin{enumerate}
\item[(i)]
From the $F \in \mathcal{L}_A$, we obtain $F_A \in \mathcal{L}$, that means, because of closure property of class $ \mathcal{L}$ with respect the convolution it holds  $F_A^{n*}\in \mathcal{L}$, see \cite[Th. 3(b)]{embrechts:goldie:1980}. Since $F^{n*}\in \mathcal{S}_A$, by Lemma \ref{lem.KP.3.a} is implied the $\overline{F_A^{n*}}(x) \asymp F^{n*}(x\,A)$. Therefore, 
$\overline{F_A'} \asymp \overline{F_A^{n*}}(x)$, $F_A' \in \mathcal{S}$ and $F_A^{n*} \in \mathcal{L}$. So we find $F^{n*}_A \in \mathcal{S}$, from \cite[Th. 2.1(a)]{kluppelberg:1988}. But we know that $\mathcal{S}$ is closed with respect to convolution roots from \cite[Th. 2]{embrechts:goldie:veraverbeke:1979}, so $F_A \in \mathcal{S}$, from where we find $F \in \mathcal{S}_A$.
\item[(ii)]
Since $F^{n*} \in \mathcal{A}_A \subsetneq \mathcal{S}_A$, $F \in \mathcal{L}_A$, then form part $(i)$ we obtain  $F \in \mathcal{S}_A$. From \cite[Cor. 4.10]{samorodnitsky:sun:2016}, we find that $F^{n*}(x\,A) \sim n\,F(x\,A)$, which together with Proposition \ref{pr.KP.2.2} implies $F \in \mathcal{A}_A $.~\halmos
\end{enumerate}

Finally, we give a result on infinite divisibility of random vectors from classes $\mathcal{S}_A,\,\mathcal{A}_A,\,\mathcal{D}_A$. Let us remind that if a distribution $F_A$ with support on $\bbr_+$ is infinitely divisible, then its Laplace transform is of the form
\beao
\widehat{F}_A(s)= \int_{0-}^\infty e^{-s\,x}\,F_A(dx) = \exp\left\{-a\,s - \int_{0}^\infty (1-e^{-s\,x})\,\nu_A(dx) \right\}\,,
\eeao
for any $s\geq 0$, where $a\geq 0$ and $\nu_A$ is the L\'{e}vy measure, with tail $\overline{\nu_A}(x) := \nu(x,\,\infty)$, such that it holds $\overline{\nu}_A(1) := \nu_A(1,\,\infty)< \infty$ and
\beao
\int_0^1 x\,\nu_A(dx) < \infty\,.
\eeao 
Let denote the normalized L\'{e}vy measure of $F_A$ as
\beao
\nu_{A,1}(x) := \dfrac{\nu_A(x)}{\overline{\nu}_A(1)}\,{\bf 1}_{\{x> 1\}}\,.
\eeao

We focus our interest on infinitely divisible distribution with heavy tails to the asymptotic behavior of the L\'{e}vy measure and its tails, see \cite{embrechts:goldie:veraverbeke:1979}, \cite{shimura:watanabe:2005} and \cite{watanabe:yamamuro:2010}. In multivariate set up we mention the papers \cite{hult:lindskog:2006} for $MRV$ distribution class.

\bpr \label{pr.KP.2.4}
Let  $A \in \mathscr{R}$ be some fixed set and ${\bf X}$ be random vector with infinitely divisible distribution $F$. 

(1) The following are equivalent
\begin{enumerate}
\item[(i)]
$F \in\mathcal{S}_A$, 
\item[(ii)]
$\nu_{A,1} \in \mathcal{S}$,
\item[(iii)]
$\bF_A(x) \sim \overline{\nu}_A(x)$.
\end{enumerate}

(2) The following are equivalent
\begin{enumerate}
\item[(i)]
$F \in\mathcal{A}_A$, 
\item[(ii)]
$\nu_{A,1} \in \mathcal{A}$.
\end{enumerate}

(3) The following are equivalent
\begin{enumerate}
\item[(i)]
$F \in \mathcal{D}_A$, 
\item[(ii)]
$\nu_{A,1} \in \mathcal{D}$,
\end{enumerate}
Further, if $\nu_{A,1} \in \mathcal{D}$, then it holds $\bF_A(x) \asymp \overline{\nu}_A(x)$.
\epr

\pr~
Since ${\bf X}$ is a non-negative random vector of infinitely divisible $d$-variate distribution, all the non-negative and non-degenerated at zero linear combinations of the components of ${\bf X}$, should follow infinitely divisible 
one-dimensional distribution, see for example in \cite[Th. 3.2]{horn:steutel:1978}. Hence, the random variable $Y_A$ follows some infinitely divisible distribution $F_A$. Hence, the assertions (1) and (2), are implied by \cite[Th. 1]{embrechts:goldie:veraverbeke:1979}, and \cite[Th. 4.1]{konstantinides:passalidis:2024d}, respectively (recall also relations \eqref{eq.KP.2.15}, \eqref{eq.KP.2.16}). The assertion (3) follows by \cite[Prop. 4.1 (i)]{watanabe:1996}, while the last statement comes from \cite[Prop. 4.1 (ii)]{watanabe:1996}, see also \cite[Th. C]{shimura:watanabe:2005}.
~\halmos

\section{Dependence scale mixtures} \label{sec.KP.4}

Now we examine the closure properties of several multivariate distribution classes with respect to dependent scale mixtures, as introduced in relation \eqref{eq.KP.1.1}. After that, we study the presence of multivariate linear single big jump in scaled mixture sums of relation \eqref{eq.KP.1.2}.

\subsection{Closure properties with repsect to scale mixtures.} \label{sec.KP.4.1}

In one-dimensional heavy-tailed distributions an important question in practical applications is the closure property with respect to product (convolution), namely if $Z_1 \stackrel{d}{\sim} V \in \mathcal{B}$ for some distribution class $\mathcal{B}$, then under what conditions the distribution of the product $Z_1 \cdot Z_2 \stackrel{d}{\sim} H$ belongs to $\mathcal{B}$ too? See for example the seminal papers \cite{cline:samorodnitsky:1994}, \cite{tang:2006}, \cite{tang:2008b}, \cite{xu:cheng:wang:cheng:2018}, \cite{cui:wang:2020} and \cite{konstantinides:leipus:siaulys:2022} for independent $Z_1$, $Z_2$, and \cite{yang:wang:2013}, \cite{chen:xu:cheng:2019}, \cite{konstantinides:passalidis:2024b} for dependent $Z_1$, $Z_2$. See also in \cite[Ch. 5]{leipus:siaulys:konstantinides:2023} for more discussions on this topic.

Let us define as scale mixture the vector 
\beam \label{eq.KP.3.1} 
\Theta\,{\bf X}=(\Theta\,X_{1},\,\ldots,\,\Theta\,X_{d})^{\top}\,,
\eeam 
where $\Theta$ represents a non-negative, non-degenerate at zero random variable. In multivariate case, for vector ${\bf X} = (X_1,\,\ldots,\,X_d)^{\top}\stackrel{d}{\sim}F$, we are interested in closure property with respect to distribution of the scale mixture $\Theta\,{\bf X}$. Scale mixtures of \eqref{eq.KP.3.1} is an important tool in risk theory and risk management applications. For example we can consider an insurer, who operates $d$-lines of 
business with corresponding random claims $X_{1},\,\ldots,\,X_{d}$ over a concrete time horizon. Further, the modern insurance companies are forced, by competition, to invest their surplus into risk-free and risky assets. Therefore the coefficient $\Theta$ plays the role of the stochastic discount factor, which, in the simplest case of constant interest rate, is degenerate to a positive value.

Many papers studied the asymptotic behavior of risk measures for the quantity
\beao
\sum_{i=1}^d \Theta\,X_i\,,
\eeao
see for example \cite{zhu:li:2012}. From \cite{li:sun:2009} we obtain a detailed description of the application of relation \eqref{eq.KP.3.1}, when the vector ${\bf X}$ follows a multivariate heavy-tailed distribution. We refer the reader to \cite[Sec. 7.3.3]{mcneil:frey:embrechts:2005} for its applications in risk management.

Thus, we study the following question: If ${\bf X} \stackrel{d}{\sim} F \in \mathcal{B}_A$, can we find sufficient conditions, in order to have $\Theta\,{\bf X} \stackrel{d}{\sim} G \in \mathcal{B}_A$? This is the closure property of class $\mathcal{B}_A$ with respect to scale mixture. This same question was examined when the vector ${\bf X}$ has distribution that belongs to class $MRV$, see for example \cite{basrak:davis:mikosch:2002b}, \cite{hult:samorodnitsky:2008}, \cite{fougeres:mercadier:2012}, where was examined the Hadamard product, which contains the scale mixture of relation \eqref{eq.KP.3.1} as special case. In \cite{konstantinides:passalidis:2024g} the closure property with respect to scale mixture was examined for the classes $\mathcal{C}_A$, $(\mathcal{D} \cap \mathcal{L})_A$, $\mathcal{S}_A$, $\mathcal{L}_A$, for any $A \in \mathscr{R}$, under the conditions of independence between $\Theta$ and ${\bf X}$ and some other general assumptions, see Assumption \ref{ass.KP.3.2} below. However, as was mentioned in  \cite{fougeres:mercadier:2012}, the independence condition is a rather strong requirement for practical orientation. For this reason we focus on closure properties under a dependence structure, that contains the independence as special case, see Assumption \ref{ass.KP.3.1} below.

For any set $A \in  \mathscr{R}$, let us define the random variable
\beam \label{eq.KP.3.2} 
M_{A}:=\sup\{ u\;:\;\Theta\,{\bf X} \in u\,A \}\,,
\eeam
with $M_A \stackrel{d}{\sim} G_A$. From relations \eqref{eq.KP.3.2}, \eqref{eq.KP.2.16}, we observe that
\beam \label{eq.KP.3.3} 
\PP[M_{A}>x]=\PP\left[\Theta\,Y_A> x\right]\,,
\eeam
for any $x>0$. Therefore, through relation \eqref{eq.KP.3.3}, the closure properties with respect to scale mixture, are reduced to the corresponding one-dimensional properties. The following dependence structure is based on conditional dependence, that was introduced in \cite{asimit:jones:2008}, via copulas and found a numerous of applications in risk theory, see for example \cite{asimit:badescu:2010}, \cite{li:tang:wu:2010}, \cite{yang:wang:leipus:siaulys:2013}  etc.

\begin{assumption} \label{ass.KP.3.1}
Let  $A \in \mathscr{R}$ and ${\bf X} \stackrel{d}{\sim} F$. We assume that the pair $({\bf X},\,\Theta)$ is conditionally dependent on $A$, symbolically $CD_A$, which means that $(Y_A,\,\Theta)$ is conditionally dependent, namely there exists some measurable function $h\;:\;[0,\,\infty) \to (0,\,\infty)$, that is bounded away from 0, such that it holds
\beam \label{eq.KP.3.4} 
\PP\left[Y_A>x\;|\; \Theta = t\right] \sim h(t)\,\PP\left[Y_A> x\right]\,,
\eeam
uniformly for $t\in S_{\Theta}$, where $S_{\Theta}$ represents the support of the distribution of $\Theta$.
\end{assumption}

\bre \label{rem.KP.3.1}
From relation \eqref{eq.KP.3.4} we obtain that Assumption \ref{ass.KP.3.1} includes the independence on $A$ for the pair $({\bf X},\,\Theta)$, as special case. When $Y_A$ is independent of $\Theta$, then $h(t)\equiv 1$ for any $t\in S_{\Theta} $, but the opposite is not always true see for example in \cite[Prop. 2.6]{cui:wang:2024}. We note that the uniformity of \eqref{eq.KP.3.4} is understood in the sense
\beao 
\lim_{\xto} \sup_{t\in S_{\Theta}}\left| \dfrac{\PP\left[Y_A>x\;|\; \Theta = t\right]}{ h(t)\,\PP\left[Y_A> x\right]} -1 \right|=0\,.
\eeao
From \eqref{eq.KP.2.16} we find out that the relation \eqref{eq.KP.3.4} is equivalent to 
\beam \label{eq.KP.3.6} 
\PP\left[{\bf X} \in x\,A\;|\; \Theta = t\right] \sim h(t)\,\PP\left[{\bf X} \in x\,A\right]\,,
\eeam
uniformly for $t\in S_{\Theta}$. If it holds  $t\notin S_{\Theta}$, then the probabilty in the left hand side of \eqref{eq.KP.3.4} (or, equivalently, of \eqref{eq.KP.3.6}), is understood as unconditional probability. The structure of \eqref{eq.KP.3.4} contains in general nested dependence structures, since are involved dependencies of $d+1$ random variables. We refer the reader to \cite[Exam. 5.3]{konstantinides:passalidis:yuan:2026}, for an example with nested Archimedean copulas. 
\ere

The next assumption is very common in product convolution of heavy-tailed distribution in one-dimensional case, see for example \cite{cline:samorodnitsky:1994} or \cite{tang:2006}.
\begin{assumption} \label{ass.KP.3.2}
Let some fixed $A \in \mathscr{R}$. We assume that for any  $c> 0$ it holds
\beao
\PP\left[\Theta>c\,x\right]= o\left(\PP\left[M_A> x\right]\right)\,.
\eeao
\end{assumption}

\bre \label{rem.KP.3.2}
From relation \eqref{eq.KP.3.3} we obtain that Assumption \ref{ass.KP.3.2} is equivalent to relation 
\beam \label{eq.KP.3.8} 
\PP[\Theta>c\,x] =o\left(\PP[\Theta\,Y_A> x] \right)\,,
\eeam
for any $c> 0$. When $Y_A$ follows heavy-tailed distribution, relation \eqref{eq.KP.3.8} represents a rather weak condition, see for example comments after \cite[Th. 2.1]{tang:2006}. It is easy to observe that if $\Theta$ has bounded from above support, then relation \eqref{eq.KP.3.8} follows immediately, when $Y_A$ has distribution with infinite right endpoint. 
\ere

In the following result, we extend some parts from \cite[Th. 3.2]{konstantinides:passalidis:2024g}, with respect to conditional dependence.

\bth \label{th.KP.3.1}
Let some fixed $A \in \mathscr{R}$. If the pair $({\bf X},\,\Theta)$ satisfies Assumptions \ref{ass.KP.3.1} and \ref{ass.KP.3.2}, with $F \in \mathcal{B}_A$, then for scale mixture $\Theta\,{\bf X} \stackrel{d}{\sim} G$, it holds $G \in \mathcal{B}_A$, where $\mathcal{B} \in \{\mathcal{C},\,\mathcal{D},\,(\mathcal{D}\cap\mathcal{L}),\,\mathcal{S},\,\mathcal{L},\,\mathcal{P_D},\,\mathcal{D}\cap\mathcal{P_D},\,\mathcal{A},\,\mathcal{T},\,\mathcal{D}\cap\mathcal{A},\,\mathcal{C}\cap\mathcal{P_D}\}$.
\ethe

\pr~
It is enough to show that $G_A \in \mathcal{B}$, which through relation \eqref{eq.KP.3.3} is equivalent to show that the distribution of $\Theta\,Y_A$ belongs to $\mathcal{B}$. Let introduce the collective statement 
\begin{enumerate}
\item[(i)]
For the class $\mathcal{C}$, the result is implied by \cite[Th. 2.5(3)]{konstantinides:passalidis:2024b}.
\item[(ii)]
For the class $\mathcal{D}$, the result is implied by \cite[Th. 2.5(1)]{konstantinides:passalidis:2024b}.
\item[(iii)]
For the class $\mathcal{D}\cap\mathcal{L}$, the result is implied by \cite[Th. 2.5(2)]{konstantinides:passalidis:2024b}.
\item[(iv)]
For the class $\mathcal{S}$, the result is implied by \cite[Th. 2.1(ii)]{chen:xu:cheng:2019}, for $\gamma=0$.
\item[(v)]
For the class $\mathcal{L}$, the result is implied by \cite[Th. 2.1(i)]{chen:xu:cheng:2019}, for $\gamma=0$.
\item[(vi)]
For the class $\mathcal{P_D}$, the result is implied by \cite[Th. 2.2]{konstantinides:passalidis:2024b}.
\item[(vii)]
For the class $\mathcal{D}\cap\mathcal{P_D}$, the result is implied by points (ii) and (vi).
\item[(viii)]
For the class $\mathcal{A}$, the result is implied by points (iv) and (vi).
\item[(ix)]
For the class $\mathcal{T}$, the result is implied by points (v) and (vi).
\item[(x)]
For the class $\mathcal{D}\cap\mathcal{A}$, the result is implied by points (ii) and (viii)
\item[(xi)]
For the class $\mathcal{C}\cap\mathcal{P_D}$, the result is implied by points (i) and (iv).
\end{enumerate}

Since $F \in \mathcal{B}_A$, we get $Y_A \stackrel{d}{\sim} F_A \in \mathcal{B}$. Thus, by Assumptions \ref{ass.KP.3.1} and \ref{ass.KP.3.2}, see also \eqref{eq.KP.3.8}, via the previous collective statement we obtain $\Theta\,Y_A \stackrel{d}{\sim} G_A \in \mathcal{B}$.   
~\halmos

If instead of $A$ we put $\mathscr{R}$, Theorem \ref{th.KP.3.1} remains valid (in the case when Assumptions \ref{ass.KP.3.1} and \ref{ass.KP.3.2} hold for any $A \in  \mathscr{R}$). 

Next, we present an extension, in some sense, of the Breiman's theorem in multivariate set up, which does not surpass or is covered by the multivariate versions in \cite{basrak:davis:mikosch:2002b}, \cite{hult:samorodnitsky:2008}, \cite{fougeres:mercadier:2012}, since from  one side is restricted to set from family $\mathscr{R}$, and from the other side provides a dependence structure between ${\bf X}$ and $\Theta$, that was not examined before. So we get a tool for more direct calculation of the asymptotic expressions, in the case $F \in MRV$, under some extra assumptions. As it seems that it is immediately connected with one-dimensional properties and concretely with Breiman's theorem, see in \cite{breiman:1965}, \cite{cline:samorodnitsky:1994}, \cite{yang:wang:2013}, \cite{cui:wang:2024} for independent or dependent versions of Breiman's lemma.

\bpr \label{pr.KP.3.1}
Let $A \in \mathscr{R}$ be some fixed set. If the pair $({\bf X},\,\Theta)$ satisfies Assumption \ref{ass.KP.3.1}, with $F \in MRV(\alpha,\,V,\,\mu)$, with $\alpha \in (0,\,\infty)$ and $\mu(A)>0$. We assume that it holds $\E\left[\Theta^p\,h(\Theta) \right] < \infty$, for some $p > \alpha$, and
\beam \label{eq.KP.3.10a} 
\PP[\Theta>x] =o\left(\PP[Y_A> x] \right)\,,
\eeam
then it holds
\beam \label{eq.KP.3.10} 
\PP[\Theta\,{\bf X} \in x\,A]\sim\E\left[\Theta^{\alpha}\,h(\Theta) \right] \,\mu(A)\,\overline{V}(x) \,,
\eeam
and furthermore $G_A \in \mathcal{R}_{-\alpha}$.
\epr

\pr~
From Assumption \ref{ass.KP.3.1}, \eqref{eq.KP.3.10a}, the moment condition $\E\left[\Theta^p\,h(\Theta) \right] < \infty$ for some $p > \alpha$, relation \eqref{eq.KP.3.3} and inclusion $F_A \in \mathcal{R}_{-\alpha}$ (see Remark \ref{rem.KP.3.A}) through the \cite[Th. 2.1(i)]{cui:wang:2024}, we obtain $M_A \stackrel{d}{\sim} G_A \in \mathcal{R}_{-\alpha}$ and it holds
\beao \label{eq.KP.3.11} 
\PP[\Theta\,Y_A > x]\sim \E\left[\Theta^{\alpha}\,h(\Theta) \right] \,\PP[Y_A> x]\sim \E\left[\Theta^{\alpha}\,h(\Theta) \right] \, \mu(A)\,\overline{V}(x) \,.
\eeao
~\halmos

Indeed, due to Remark \ref{rem.KP.3.A}, from Proposition \ref{pr.KP.3.1} does NOT follow $G \in MRV(\alpha,\,V,\,\mu)$. We note also that due to condition $\mu(A) >0$ (and since by definition $\mu(A) < \infty$, for any $A \in \mathscr{R}$) relation \eqref{eq.KP.3.10} is equivalent to $\PP[\Theta>x] =o\left(\bV(x) \right)$.

\subsection{Scale mixture sums.} \label{sec.KP.4.2}

Now we look for the asymptotic behavior of the scale mixture sums from \eqref{eq.KP.1.2} and concretely we are interested in multivariate linear single big jump on set $xA$, namely
\beam \label{eq.KP.3.11b} 
\PP\left[{\bf S}_n^{\Theta} \in x\,A \right] \sim \sum_{i=1}^n \PP\left[\Theta^{(i)}\,{\bf X}^{(i)} \in x\,A \right] \,.
\eeam
Let us notice that the sets $A$ can be connected with some ruin-sets, see Sec. \ref{sec.KP.5}. Relation \eqref{eq.KP.3.11b} was studied in non-weighted form, namely when $\Theta^{(i)}\equiv 1$ for $i=1,\,\ldots,\,n$, in several papers under independent and identically distributed random vectors ${\bf X}^{(i)}$ with distribution from $MRV$. In \cite[Sect. 4]{das:fashenhartmann:2023} there is extension of such results in case of not necessarily identically distributed ${\bf X}^{(i)}$ and in \cite[Sect. 4.1]{konstantinides:passalidis:2024g} the classes $\mathcal{C}_A$, $(\mathcal{D}\cap \mathcal{L})_A$, $\mathcal{S}_A$, with corresponding random vectors satisfing some weak dependence structures, can be found. We note that, all the random vectors contain arbitrarily dependent components. 

Now, we consider that the pairs $({\bf X}^{(1)},\,\Theta^{(1)}),\,\ldots,\,({\bf X}^{(n)},\,\Theta^{(n)})$, are independent and each pair $({\bf X}^{(i)},\,\Theta^{(i)})$, for $i=1,\,\ldots,\,n$, $n \in \bbn$, contains the dependence structure from Assumption \ref{ass.KP.3.1} (note that each  ${\bf X}^{(i)}$ has arbitrarily dependent components). We have to accept that in case $d=1$ and $A=(1,\,\infty)$, relation \eqref{eq.KP.3.11b} represents a well-studied asymptotic formula, see for example \cite{gao:wang:2010}, \cite{cheng:2014}, \cite{yu:wang:cheng:2015}, and \cite{chen:cheng:2024}, with several weights and different dependencies.

In the following result the function $h_i$ corresponds to the similar one in Assumption \ref{ass.KP.3.1}  for the pairs $({\bf X}^{(i)},\,\Theta^{(i)})$, for $i=1,\,\ldots,\,n$ respectively. Additionally, Assumption \ref{ass.KP.3.2} should be understood as 
\beao
\PP\left[\Theta^{(i)} > c\,x \right] =o\left(\PP\left[M_A^{(i)}> x\right] \right)\,,
\eeao 
for any $c>0$, with $M_A^{(i)}:=\sup\left\{ u\;:\;\Theta^{(i)}\,{\bf X}^{(i)} \in u\,A \right\}$, for $i=1,\,\ldots,\,n$.

\bth \label{th.KP.3.2} 
Let $A \in \mathscr{R}$ be some fixed set. We suppose that the sequence of pairs $({\bf X}^{(1)},\,\Theta^{(1)}),\,\ldots,\,({\bf X}^{(n)},\,\Theta^{(n)})$, $n \in \bbn$, represent independent random vectors, and each pair $({\bf X}^{(i)},\,\Theta^{(i)})$, satisfies Assumption \ref{ass.KP.3.1}, for $i=1,\,\ldots,\,n$. In this case, we obtain
\begin{enumerate}
\item[(i)]
If ${\bf X}^{(1)},\,\ldots,\,{\bf X}^{(n)}$ follow a common distribution $F \in \mathcal{S}_A$ and the random variables $\Theta^{(i)}$ are such that $\PP[\Theta^{(i)} > 0]> 0$, and $\PP\left[\Theta^{(i)} \in [0,\,b]\right] =1$, for some $b\in (0,\,\infty)$ and for any $i=1,\,\ldots,\,n$. Then \eqref{eq.KP.3.11b} holds. 
\item[(ii)]
If ${\bf X}^{(1)},\,\ldots,\,{\bf X}^{(n)}$ follow the distributions $F_{1},\,\ldots,\,F_{n} \in (\mathcal{D} \cap \mathcal{L})_A$ respectively, and for each pair $({\bf X}^{(i)},\,\Theta^{(i)})$ Assumption \ref{ass.KP.3.2} is true. Then \eqref{eq.KP.3.11b} holds. 
\end{enumerate}
\ethe

\pr~
For the upper bound of \eqref{eq.KP.3.11b} in both cases we obtain
\beam \label{eq.KP.3.12} \notag
&&\PP[{\bf S}_n^{\Theta} \in x\,A ] = \PP\left[ \sup_{{\bf p}\in I_A}  {\bf p}^T\,\left( \sum_{i=1}^n \Theta^{(i)} {\bf X}^{(i)} \right) >x\right]\leq \PP\left[\sum_{i=1}^n \left(  \sup_{{\bf p}\in I_A}  {\bf p}^T \Theta^{(i)} {\bf X}^{(i)} \right) >x\right]\\[2mm] 
&&= \PP\left[\sum_{i=1}^n \Theta^{(i)}\,Y_A^{(i)} > x \right] \sim \sum_{i=1}^n \PP\left[\Theta^{(i)}\,Y_A^{(i)} > x \right] \sim  \sum_{i=1}^n \PP\left[\Theta^{(i)}\,{\bf X}^{(i)} \in x\,A \right]\,,
\eeam
where in the pre-last step, for the case $(i)$, we used \cite[Th. 1.2]{yang:leipus:siaulys:2012}, while for the case $(ii)$, we needed \cite[Exam. 4.1(2)]{konstantinides:passalidis:2024b}, that generalizes \cite[Th. 1]{yang:wang:leipus:siaulys:2013}.

For the lower bound, in both cases, taking into account that the $\Theta^{(i)}$  and ${\bf X}^{(i)}$ are non-negative and the $xA$ is 'increasing set', by Bonferroni's inequality we find out that
\beam \label{eq.KP.3.13} \notag
&&\PP\left[{\bf S}_n^{\Theta} \in x\,A \right] \geq \PP\left[ \bigcup_{i=1}^n  \left\{\Theta^{(i)}\,{\bf X}^{(i)} \in x\,A \right\} \right]\\[2mm] \notag
&& \geq \sum_{i=1}^n\PP\left[ \Theta^{(i)}\,{\bf X}^{(i)} \in x\,A \right]-\sum_{1\leq i<j \leq n} \PP\left[ \Theta^{(i)}\,{\bf X}^{(i)} \in x\,A\,,\;\Theta^{(j)}\,{\bf X}^{(j)} \in x\,A\right]\\[2mm]
&&= \sum_{i=1}^n\PP\left[ \Theta^{(i)}\,Y_A^{(i)} > x \right] - \sum_{1\leq i<j \leq n} \PP\left[\Theta^{(i)}\,Y_A^{(i)} > x \right] \, \PP\left[\Theta^{(j)}\,Y_A^{(j)} > x \right] \\[2mm] \notag
&&\sim  \sum_{i=1}^n\PP\left[ \Theta^{(i)}\,Y_A^{(i)} > x \right] \sim  \sum_{i=1}^n \PP\left[\Theta^{(i)}\,{\bf X}^{(i)} \in x\,A \right]\,,
\eeam
where in the third step we used relations \eqref{eq.KP.3.2} and \eqref{eq.KP.3.3} and the fact that the pairs $(Y_A^{(i)},\,\Theta^{(i)})$, $(Y_A^{(j)},\,\Theta^{(j)})$ are independent for any $1\leq i \neq j \leq n$. From relations \eqref{eq.KP.3.12} and \eqref{eq.KP.3.13} we conclude \eqref{eq.KP.3.11b}.
~\halmos

It is worth to observe that while in the first point of Theorem \ref{th.KP.3.2} we use subexponential distribution class $\mathcal{S}_A$, and simultaneously we assume identically distributed ${\bf X}^{(i)}$, and the $\Theta^{(i)}$ are bounded from above, in the second point we reduce the distribution class to $(\mathcal{D} \cap \mathcal{L})_A$ but we permit non identically distributed ${\bf X}^{(i)}$, and more general conditions for $\Theta^{(i)}$. Next, a further restriction of the distribution class to $MRV$, provides a more direct form of relation \eqref{eq.KP.3.11b}.
   
\bco \label{cor.KP.3.1}
Let $A \in \mathscr{R}$ be some fixed set. We suppose that the conditions of Theorem \ref{th.KP.3.2}(ii) are valid under the restriction that $F_i \in MRV(\alpha_i,\,V_i,\,\mu_i)$, with $\alpha_i \in (0,\,\infty)$, $\mu_i(A) >0$, and there are some $p_i > \alpha_i$, such that $\E\left[(\Theta^{(i)})^{p_i}\,h_i(\Theta^{(i)}) \right] < \infty$, for any $i=1,\,\ldots,\,n$. Also we assume that each pair $({\bf X}^{(i)},\,\Theta^{(i)})$, for $i=1,\,\ldots,\,n$, satisfies the relation \eqref{eq.KP.3.10a}.  Then
\beao
\PP\left[{\bf S}_n^{\Theta} \in x\,A \right] \sim \sum_{i=1}^n \E\left[\left(\Theta^{(i)}\right)^{\alpha_i}\,h_i\left(\Theta^{(i)}\right) \right]\, \mu_i(A)\,\overline{V}_i(x)\,.
\eeao
\eco

\pr~
We use Theorem \ref{th.KP.3.2}(2), together with the fact that for $\mu_i(A) >0$ it holds $MRV(\alpha_i,\,V_i,\,\mu_i) \subsetneq (\mathcal{D} \cap \mathcal{L})_A$, and through Proposition \ref{pr.KP.3.1}.
~\halmos

\bre \label{rem.KP.3.4}
According to \cite[Prop. 2.4]{cui:wang:2024}, function $h$ is bounded from above, that means there exists $K>0$, such that it holds
\beam \label{eq.KP.3.15}
h(t) \leq K\,,
\eeam
for any $t \in S_{\Theta}$. By \eqref{eq.KP.3.15} we find that if $\E\left[\Theta^{p}\right] < \infty$ for some $p > \alpha$, then 
\beao
\E\left[\Theta^{p}\,h(\Theta) \right] =\int_{S_{\Theta}} t^p\,h(t)\,\PP[\Theta \in dt ] \leq K\,\E\left[\Theta^{p}\right] < \infty\,.
\eeao
Hence, the moment conditions in Proposition \ref{pr.KP.3.1} and Corollary \ref{cor.KP.3.1} can be checked directly. However, $\E\left[\Theta^{p}\,h(\Theta) \right] < \infty $ does not imply necessarily $\E\left[\Theta^{p}\right] < \infty$, see \cite[Exam. 2.1]{cui:wang:2024}.
\ere

\section{Finite-time ruin probability in a time-dependent risk model} \label{sec.KP.5}

Recently, because of increasing popularity of dependence in the frame of insurance industry, the number of papers on multivariate risk models constantly increases. The bivariate models are more attractive, since they present the less mathematical complexity, than multivariate one, see for example \cite{chen:wang:wang:2013}, \cite{yang:li:2014}, \cite{jiang:wang:chen:xu:2015}, \cite{yang:chen:yuen:2024}, among others. From the other side, the multivariate models are restricted in the $MRV$ condition for the claims, and even more usually each claim vector has asymptotically dependent components, namely the Radon measure is such that it holds $\mu\left((1,\,\infty]\times \cdots \times (1,\,\infty] \right)>0$, see for example \cite{yang:su:2023}, \cite{cheng:konstantinides:wang:2024}. In larger than $MRV$ distribution classes on claims, we find \cite{samorodnitsky:sun:2016}, \cite{konstantinides:passalidis:2024g}, \cite{konstantinides:liu:passalidis:2026}, where the multivariate subexponential distributions were studied. In one-dimensional risk models, after 2010, the time-dependent risk models are well studied, see for example \cite{asimit:badescu:2010}, \cite{li:tang:wu:2010} for some seminal papers on this topic. In multivariate set up we know only \cite{li:2016} and \cite{cheng:konstantinides:wang:2022}, for time-dependent risk models, with claim-vector distributions from class $MRV$ and asymptotically dependent components.

Next, we examine the multivariate risk model with constant interest force, whose arrival counting process represents a common Poisson process, and it permits a dependence structure among the multivariate claims and among the inter-arrival times.
More precisely, we will study the following model:
\beam \label{eq.KP.4.18}
{\bf U}(t) := x\,\left( 
\begin{array}{c}
	l_{1} \\ 
	\vdots \\ 
	l_{d} 
\end{array} 
\right)
+ \left( 
\begin{array}{c}
	\int_{0}^t e^{-r\,y}\,c_1(y)\,dy \\ 
	\vdots \\ 
	\int_{0}^t e^{-r\,y}\,c_d(y)\,dy 
\end{array} 
\right)-{\bf D}_r(t)\,,
\eeam
for any $t\geq 0$, where $x$ represents the total initial capital, and the $l_1,\,\ldots,\,l_d >0$ are weights with 
\beao
\sum_{i=1}^d l_i =1\,,
\eeao 
and they depict the allocation of the initial capital to $d$-lines of business. Further, the $c_i(s)$, $i=1,\,\ldots,\,d $, are premium density of the $i$-th line of business, for which there exists a constant $M_i>0$, such that it holds $c_i(s) \in [0,\,M_i]$, for any $s \in [0,\,t]$, with $c_i(0)=0$. This mild boundness condition on the premium-densities, allows us are arbitrary dependence between $c_i(t)$ ($i=1,\,\ldots,\,d$) and ${\bf D}_r(t)$.

${\bf D}_r(t)$ depicts the discounted aggregate claims from the side of insurer, described through the following relation:
\beam \label{eq.KP.4.1}
{\bf D}_r(t) := \sum_{i=1}^{N(t)} {\bf X}^{(i)}\,e^{-r\,\tau_i}= \left( \sum_{i=1}^{N(t)} X_{1}^{(i)}\,e^{-r\,\tau_{i} } ,\,\ldots ,\,\sum_{i=1}^{N(t)} X_{d}^{(i)}\,e^{-r\,\tau_{i} } \right)^{\top}\,,
\eeam
for $t\geq 0$. On the discounted aggregate claims from \eqref{eq.KP.4.1}, the claim vectors 
\beao
{\bf X}^{(i)}=(X_{1}^{(i)},\,\ldots,\,X_{d}^{(i)})^{\top}\,,
\eeao 
with $i \in \bbn$, are independent and identically distributed copies of the vector  
\beao
{\bf X}=(X_{1},\,\ldots,\,X_{d})^{\top}\,,
\eeao 
 with $X_{j}^{(i)}$,depicting the $i$-th claim of the $j$-th portfolio, for any $j=1,\,\ldots,\,d$. We should mention that the vector $ {\bf X}^{(i)}$ can have up to $d-1$ zero components.

The arrival of ${\bf X}^{(i)}$ happens at time moment $\tau_i$, with the sequence $\{\tau_i\,,\;i \in \bbn\}$ to constitute a Poisson process with counting process $\{ N(t)\,,\; t\geq 0 \}$. Further, we denote by $\{ \Theta_i\,,\; i \in \bbn \}$ the sequence of the inter-arrival times, namely 
\beao
\Theta_i = \tau_i-\tau_{i-1}\,,
\eeao 
for $i \in \bbn$ and $\tau_0:=0$. The constant interest rate is denoted by $r\geq 0$. In such a model we consider that the pairs $({\bf X}^{(i)},\,\Theta_i)$ satisfy an analogue of Assumption \ref{ass.KP.3.1} as a dependence structure, see Assumption \ref{ass.KP.4.1} below, which contains the independence as special case. Namely, in each vector ${\bf X}^{(i)}$ there are arbitrarily dependent components (under the assumption that satisfy some multivariate class), and there exist a weak dependence structure among ${\bf X}^{(i)}$ and the inter-arrival times $\Theta_i\,,\; i \in \bbn$. The common counting process is also another source of dependence between the d-lines of business.

Let us firstly consider the asymptotic behavior of $\PP\left[{\bf D}_r(t) \in x\,A\right]$, for any finite $t \in [0,\,\infty)$. The set $A$ can take many interesting forms, as for example in case $d=1$ it can be 
\beao
A=(1,\,\infty)\,,
\eeao 
or in case $d\geq 2$, it can be 
\beao
A=\left\{{\bf x}\;:\; \sum_{i=1}^d c_i\,x_i> u\right\}\,,
\eeao
where $u>0$, with $c_1+\cdots+c_d =1$, that corresponds to the event when the sum of discounted aggregate claims of d-lines exceeds the initial capital of the company. In \cite{samorodnitsky:sun:2016} we find that such sets $A$ belong to the family of sets $\mathscr{R}$ and in \cite[Rem. 2.2]{konstantinides:passalidis:2024g} for possible exceptions.

Before formulating the exact assumptions of our model, we need an auxiliary result, related to scale mixture sums. In opposite to subsection \ref{sec.KP.4.2}, now $\Theta^{(i)}$ and ${\bf X}^{(i)}$ are independent.

\ble \label{lem.KP.4.1}
Let $A \in \mathscr{R}$ be some fixed set. We consider that the $ {\bf X}^{(1)},\,\ldots,\, {\bf X}^{(n)}$ are independent random vectors, with corresponding distributions $F_{1},\,\ldots,\, F_{n} \in \mathcal{L}_A$ and the $\Theta^{(1)},\,\ldots,\, \Theta^{(n)}$ non-negative, non-degenerate at zero, random variables, that are independent of $ {\bf X}^{(1)},\,\ldots,\, {\bf X}^{(n)}$. We assume that there exist some d-variate distribution $F$ such that 
\beao
F_i(x\,A) \asymp F(x\,A)\,,
\eeao 
for any $i=1,\,\ldots,\,n $.
\begin{enumerate}
\item[(i)]
If $F \in \mathcal{S}_A$ and 
\beao
\PP\left[\Theta^{(i)} \in [0,\,b_i]\right]=1\,,
\eeao 
for any $i=1,\,\ldots,\,n $, with $b_i \in (0,\,\infty)$, then \eqref{eq.KP.3.11b} holds. 
\item[(ii)]
If $F \in \mathcal{A}_A$ and Assumption \ref{ass.KP.3.2} is satisfied for any $i=1,\,\ldots,\,n $, then \eqref{eq.KP.3.11b} holds. 
\end{enumerate}
\ele

\pr~
For the upper bound of \eqref{eq.KP.3.11b}, we use relation \eqref{eq.KP.3.12}, where in pre-last step we apply \cite[Th. 1]{tang:yuan:2014} for the case $(i)$, and  \cite[Th. 2]{tang:yuan:2014} for the case $(ii)$. For the lower bound of  \eqref{eq.KP.3.11b}, we use the fact that $\Theta^{(i)}$ and ${\bf X}^{(i)}$ are non-negative, and the $xA$ is increasing set. Hence, we obtain
\beao
&&\PP\left[{\bf S}_n^{\Theta} \in x\,A \right] \geq \PP\left[ \bigcup_{i=1}^n  \left\{\Theta^{(i)}\,{\bf X}^{(i)} \in x\,A \right\} \right] \\[2mm]
&&\geq \sum_{i=1}^n \PP\left[ \Theta^{(i)}\,{\bf X}^{(i)} \in x\,A \right]-\sum_{1\leq i<j \leq n} \PP\left[ \Theta^{(i)}\,{\bf X}^{(i)} \in x\,A\,,\;\Theta^{(j)}\,{\bf X}^{(j)} \in x\,A\right]\\[2mm] 
&&= \sum_{i=1}^n \PP\left[ \Theta^{(i)}\,Y_A^{(i)} > x \right] - \sum_{1\leq i<j \leq n} \PP\left[ \Theta^{(i)}\,{\bf X}^{(i)} \in x\,A\,,\;\Theta^{(j)}\,{\bf X}^{(j)} \in x\,A\right] \\[2mm]  
&&\sim  \sum_{i=1}^n\PP\left[ \Theta^{(i)}\,Y_A^{(i)} > x \right] \sim  \sum_{i=1}^n \PP\left[\Theta^{(i)}\,{\bf X}^{(i)} \in x\,A \right]\,,
\eeao
where in the pre-last step we took into account relation \cite[eq. (15)]{tang:yuan:2014} for both points.
~\halmos

Next, we present a crucial Assumption for our model. As we can see Assumption \ref{ass.KP.4.1} is used to model the dependence between vector ${\bf X}^{(i)}$ and inter-arrival time $\Theta_i$.

\begin{assumption} \label{ass.KP.4.1}
Let $A \in \mathscr{R}$ be some fixed set. We assume that $({\bf X}^{(i)},\,\Theta_i)$ are independent and identically distributed random vectors, for  any $i\in \bbn$. Further, we suppose that there exists a positive, measurable function $h\;:\;[0,\,\infty) \to (0,\,\infty)$, such that it holds
\beam \label{eq.KP.4.3} 
\PP\left[ {\bf X} \in x\,A\;|\;\Theta= s \right]\sim h(s)\,\PP\left[  {\bf X} \in x\,A\right]\,,
\eeam
uniformly for $s \in (0,\,t]$, for any fixed $t>0$, with $\inf_{s\in(0,\,t]}h(s)>0$, conventionally. Here $({\bf X},\,\Theta)$ is the generic random vector of $\left({\bf X}^{(i)},\,\Theta_i\right)$, $i\in\bbn$.
\end{assumption}

In what follows we define 
\beao
\Lambda_n(t):=\left\{ {\bf s} \in (0,\,t)^n\;:\; \sum_{i=1}^n s_i < t \right\}\,.
\eeao 
Since the counting process $\{N(t)\,,\;t\geq 0\}$ represents a Poisson process, we know that it holds
\beam \label{eq.KP.4.4} 
\PP\left[ {\bf \Theta}={\bf s}\;|\;N(t)= n \right]=\dfrac{n!}{t^n}\,,
\eeam
for all ${\bf s}\in\Lambda_n(t)$, see for example \cite[p. 188]{embrechts:kluepellberg:mikosch:1997}. The next theorem generalize the main results in \cite{asimit:badescu:2010} in multivariate setup. Let us observe that the theorem is identical with \cite[Th. 3.1, 3.2]{asimit:badescu:2010}, for $d=1$ and $A=(1,\,\infty)$. We also observe that in opposite to \cite[Theorem 2.1]{li:2016}, where was studied a time-dependent renewal risk model, here we do not restricted only on asymptotic dependence assumption among the components of the claim vectors, while we are restricted to Poisson risk model and constant interest force.

\bth \label{th.KP.4.1}
Let $A \in \mathscr{R}$ be some fixed set. We consider the compound Poisson process \eqref{eq.KP.4.1}, under the Assumption \ref{ass.KP.4.1}, for any $s \in (0,\,t]$, where $t$ is a fixed number from interval $(0,\,\infty)$. Suppose that all ${\{\bf X}^{(i)},\, i\in\bbn\}$ follow a common distribution $F$. 
\begin{enumerate}
\item[(i)]
If $r=0$, and $F \in \mathcal{S}_A$, then it holds
\beam \label{eq.KP.4.5} 
\PP\left[ {\bf D}_0(t) \in x\,A \right] \sim C_0\,\PP\left[ {\bf X} \in x\,A \right]\,,
\eeam
where
\beao
C_0=\sum_{n=1}^{\infty} \lambda^n\,e^{-\lambda\,t} \int_{\Lambda_n(t)} \sum_{i=1}^n h(s_i)\,d{\bf s}\,,
\eeao
\item[(ii)]
If $r>0$, and $F \in MRV(\alpha,\, V,\,\mu)$ for some $\alpha \in (0,\,\infty)$ and $\mu(A) >0$, then it holds
\beam \label{eq.KP.4.7}  
\PP\left[ {\bf D}_r(t) \in x\,A \right] \sim C_r\,\,\mu(A)\,\overline{V}(x)\,,
\eeam
where
\beao
C_r=\sum_{n=1}^{\infty} \lambda^n\,e^{-\lambda\,t} \int_{\Lambda_n(t)} \sum_{i=1}^n h(s_i)\,e^{-\alpha r \sum_{j=1}^i\,s_j}\,d{\bf s}\,,
\eeao
\end{enumerate}
\ethe

\pr~
At first we consider non-negative interest rate $r\geq 0$. By relation \eqref{eq.KP.4.4}, via the law of total probability, we obtain
\beam \label{eq.KP.4.8}  \notag
&&\PP\left[ {\bf D}_r(t) \in x\,A \right]= \PP\left[ \sum_{i=1}^{N(t)} {\bf X}^{(i)}\,e^{-r\,\tau_i} \in x\,A \right] \\[2mm] \notag
&&=  \PP\left[ \sum_{i=1}^{N(t)} {\bf X}^{(i)}\,\exp\left\{-r\,\sum_{j=1}^i \Theta_j \right\} \in x\,A \right] =\\[2mm] \notag
&&\sum_{n=1}^{\infty} \int_{\Lambda_n(t)} \PP\left[ \sum_{i=1}^{n} {\bf X}^{(i)} \exp\left\{-r \sum_{j=1}^i s_j \right\} \in x A \,\Big|\,{\bf \Theta}={\bf s},\;N(t)=n\right]\,\PP\left[{\bf \Theta}={\bf s},\,N(t)=n\right] d{\bf s}\\[2mm] \notag
&&=\sum_{n=1}^{\infty} \int_{\Lambda_n(t)} \PP\left[ \sum_{i=1}^{n} {\bf X}^{(i)} e^{-r \sum_{j=1}^i s_j } \in x A \,\Big|\,{\bf \Theta}={\bf s}\right]\,\PP\left[{\bf \Theta}={\bf s}\,|\,N(t)=n\right]\,\PP[N(t)=n]\,d{\bf s}\\[2mm]
&&=\sum_{n=1}^{\infty} e^{-\lambda\,t}\,\dfrac{(\lambda\,t)^n}{n!} \int_{\Lambda_n(t)} \PP\left[ \sum_{i=1}^{n} {\bf X}^{(i)} e^{-r \sum_{j=1}^i s_j } \in x A \,\Big|\,{\bf \Theta}={\bf s}\right] \,\dfrac{n!}{t^n}\,d{\bf s}\\[2mm] \notag
&&=\sum_{n=1}^{\infty} e^{-\lambda\,t}\,\lambda^n\,\int_{\Lambda_n(t)} \PP\left[ \sum_{i=1}^{n} {\bf X}^{(i)} e^{-r \sum_{j=1}^i s_j } \in x A \,\Big|\,{\bf \Theta}={\bf s}\right] \,d{\bf s}\,.
\eeam
From relation \eqref{eq.KP.4.3} and the fact that function $h$ is bounded from above, see relation \eqref{eq.KP.3.15}, we obtain that for any $x>0$ and $s_i \in (0,\,t)$, for $i\in \bbn$, there exists a constant $K^{*}\geq K \vee 1$, such that it holds
\beam \label{eq.KP.4.9} 
\PP\left[{\bf X}^{(i)}\,e^{-r \sum_{j=1}^i s_j } \in x A \,\Big|\, \Theta_{i}=s_{i}\right] &=& \PP\left[ Y_A^{(i)}\, e^{-r \sum_{j=1}^i s_j } >x\,\Big|\, \Theta_{i}=s_{i}\right]\\[2mm]  \notag
&\leq& \PP\left[ Y_A^{(i)} >x\,\Big|\, \Theta_{i}=s_{i}\right] \leq K^{*}\, \PP\left[ Y_A^{(i)} >x\right] \,.
\eeam
By relation \eqref{eq.KP.4.9} and via \cite[Lem. 3.1]{asimit:badescu:2010}, for any $\delta > 0$, there exists a positive constant $L_{\delta}>0$, such that for any $x>0$, and any $n \in \bbn$, it holds
\beam \label{eq.KP.4.10}  \notag
&&\PP\left[ \sum_{i=1}^{n} {\bf X}^{(i)} \exp\left\{-r \sum_{j=1}^i s_j \right\} \in x A \,\Big|\,{\bf \Theta}={\bf s}\right] \\[2mm]
&&=\PP\left[\sup_{{\bf p} \in I_A} {\bf p}^T\,\left(e^{-r \,s_1 }\,{\bf X}^{(1)}+\cdots+e^{-r \sum_{j=1}^n s_j }\,{\bf X}^{(n)}\right)> x \,\Big|\,{\bf \Theta}={\bf s}\right] \\[2mm] \notag
&&\leq \PP\left[ \sum_{i=1}^{n} Y_A^{(i)}\, e^{-r \sum_{j=1}^i s_j } >x\,\Big|\, {\bf \Theta}={\bf s}\right] \leq \PP\left[\sum_{i=1}^{n}  Y_A^{(i)} >x\,\Big|\,  {\bf \Theta}={\bf s}\right] \\[2mm] \notag
&&\leq L_{\delta}\,(1+\delta)^n\,\PP\left[ Y_A>x \right] =  L_{\delta}\,(1+\delta)^n\,\PP\left[ {\bf X} \in x A\right] \,.
\eeam
Therefore, from relation \eqref{eq.KP.4.10}, we can apply the dominated convergence theorem on \eqref{eq.KP.4.8}, after dividing by $\PP\left[  {\bf X} \in x A\right]$, since we can find a small enough $\vep >0$, such that it holds $\E\left[ (1+\vep)^{N(t)}\right] < \infty$, see \cite{stein:1946}.

\begin{enumerate}
\item[(i)]
Next, we consider the case $r=0$. By relation \eqref{eq.KP.4.3} and the fact that ${\bf X}^{(i)} \stackrel{d}{\sim} F \in \mathcal{S}_A$, due to the closure property of $\mathcal{S}_A$, with respect to strong 
equivalence, see \cite[Prop. 4.12(a)]{samorodnitsky:sun:2016}, we find that the conditional random vectors $\left({\bf X}^{(i)} \,\Big|\, \Theta_{i}=s_{i} \right)$, for $i\in \bbn$, are independent random vectors with distributions from class $\mathcal{S}_A$. Even more, for any $i\neq k \in \bbn$, it holds
\beao
\PP\left[  {\bf X}^{(i)}  \in x A \,\Big|\,\Theta_{i}=s_{i}\right] \asymp \PP\left[  {\bf X}^{(k)}  \in x A \,\Big|\,\Theta_{k}=s_{k}\right] \,,
\eeao
because of the fact that the pairs $({\bf X}^{(i)},\,{\bf X}^{(k)})$ are identically distributed and the pairs $(\Theta_i,\,\Theta_k)$ are also identically distributed and with the help of relation \eqref{eq.KP.4.3} (remind that $h$ is bounded away from $0$). Hence, applying Lemma \ref{lem.KP.4.1}(i) on the conditional random vectors $\left\{{\bf X}^{(i)} \,\big|\, \Theta_{i}=s_{i} \right\}$, for $i\in \bbn$, for weights degenerate to unity, via relation \eqref{eq.KP.4.8}, we obtain
\beao
&&\lim \dfrac{\PP\left[ {\bf D}_0(t) \in x\,A \right]}{\PP\left[  {\bf X} \in x A\right]} =\sum_{n=1}^{\infty} \lambda^n\,e^{-\lambda\,t}\,\int_{\Lambda_n(t)} \lim \dfrac{\PP\left[ \sum_{i=1}^{n} {\bf X}^{(i)} \in x A \,\big|\,{\bf \Theta}={\bf s}\right]}{\PP\left[  {\bf X} \in x A\right]} \,d{\bf s}= \\[2mm]
&&\sum_{n=1}^{\infty} \lambda^n\,e^{-\lambda\,t}\,\int_{\Lambda_n(t)} \lim \dfrac{ \sum_{i=1}^{n} \PP\left[{\bf X}^{(i)} \in x A \,\big|\,\Theta_{i}=s_{i}\right]}{\PP\left[  {\bf X} \in x A\right]} d{\bf s} \sim \sum_{n=1}^{\infty} \lambda^n e^{-\lambda\,t} \int_{\Lambda_n(t)} \sum_{i=1}^{n} h(s_i) d{\bf s},
\eeao
that provides relation \eqref{eq.KP.4.5}.
\item[(ii)]
Now we assume $r>0$. So we get
\beam \label{eq.KP.4.12} 
0< e^{-r\,t}\leq \exp\left\{-r \sum_{j=1}^i s_j \right\} \leq 1\,,
\eeam
and the distribution of ${\bf X}^{(i)}$ belongs to class $\mathcal{S}_A$, because $\mu(A)>0$ so $MRV \subsetneq \mathcal{S}_A $, see relation \eqref{eq.KP.2.19}. Therefore by 
Theorem \ref{th.KP.3.1}, we find that the distributions of ${\bf X}^{(i)} \exp\left\{-r \sum_{j=1}^i s_j \right\} $ also  belong to class $\mathcal{S}_A$. From Assumption \ref{ass.KP.4.1} and relation 
\eqref{eq.KP.4.12} we obtain that it holds
\beam \label{eq.KP.4.13} 
\PP\left[ {\bf X}^{(i)} \exp\left\{-r \sum_{j=1}^i s_j \right\}  \in x A \,\Bigg|\,\Theta_{i}=s_{i}\right] \sim h(s_{i}) \PP\left[  {\bf X}^{(i)} \exp\left\{-r \sum_{j=1}^i s_j \right\} \in x A\right],
\eeam
uniformly for $ s_i \in (0,\,t]$. By the closure properties of class $\mathcal{S}_A$, with respect to strong equivalence, see \cite[Prop. 4.12(a)]{samorodnitsky:sun:2016}, we find that the distributions of conditional random vectors $\left\{ {\bf X}^{(i)}\, \exp\left\{-r \sum_{j=1}^i s_j \right\}\;\Big|\;\Theta_{i}=s_{i} \right\}$ belong to class $\mathcal{S}_A$. Furthermore, by relation \eqref{eq.KP.4.13} and the fact that the vectors ${\bf X}^{(i)}$ are identically distributed with distribution from class $\mathcal{D}_A$, since $\mu(A)>0$, so $MRV \subsetneq \mathcal{D}_A $, we obtain
\beao
\PP\left[ {\bf X}^{(i)}\, \exp\left\{-r \sum_{j=1}^i s_j \right\} \in x A\;\bigg|\;\Theta_{i}=s_{i}\right] \asymp \PP\left[ {\bf X}^{(k)}\, \exp\left\{-r \sum_{j=1}^k s_j \right\} \in x A \,\bigg|\,\Theta_{k}=s_{k}\right] .
\eeao
Thus, through Lemma \ref{lem.KP.4.1}(1) and relation \eqref{eq.KP.4.8}, we can apply the dominated convergence theorem to find
\beam \label{eq.KP.4.14} \notag
&&\lim \dfrac{\PP\left[ {\bf D}_r(t) \in x\,A \right]}{\PP\left[  {\bf X} \in x A\right]} =\sum_{n=1}^{\infty} \lambda^n\,e^{-\lambda\,t}\,\int_{\Lambda_n(t)} \lim \dfrac{\PP\left[ \sum_{i=1}^{n} {\bf X}^{(i)}\,e^{-r \sum_{j=1}^i s_j } \in x A \,\big|\,{\bf \Theta}={\bf s}\right]}{\PP\left[  {\bf X} \in x A\right]} \,d{\bf s} \\[2mm] 
&&=\sum_{n=1}^{\infty} \lambda^n\,e^{-\lambda\,t}\,\int_{\Lambda_n(t)} \lim \dfrac{ \sum_{i=1}^{n} \PP\left[{\bf X}^{(i)}\,e^{-r \sum_{j=1}^i s_j } \in x A \,\big|\,{\bf \Theta}={\bf s}\right]}{\PP\left[  {\bf X} \in x A\right]} d{\bf s} \\[2mm] \notag
&&= \sum_{n=1}^{\infty} \lambda^n\,e^{-\lambda\,t}\,\int_{\Lambda_n(t)} \lim \dfrac{ \sum_{i=1}^{n} \PP\left[{\bf X}^{(i)}\,e^{-r \sum_{j=1}^i s_j } \in x A \right]\,h(s_i)}{\PP\left[  {\bf X} \in x A\right]} d{\bf s}  \\[2mm] \notag
&&= \sum_{n=1}^{\infty} \lambda^n\,e^{-\lambda\,t}\,\int_{\Lambda_n(t)} \lim \dfrac{ \sum_{i=1}^{n} h(s_i)\,\PP\left[{\bf X}^{(i)} \in x A \right]\,e^{-\alpha r \sum_{j=1}^i\,s_j}}{\PP\left[  {\bf X} \in x A\right]} d{\bf s} \\[2mm] \notag
&&= \sum_{n=1}^{\infty} \lambda^n\,e^{-\lambda\,t} \int_{\Lambda_n(t)} \sum_{i=1}^{n} h(s_i)\,e^{-\alpha r \sum_{j=1}^i\,s_j} \,d{\bf s} =C_r\,,
\eeam
where in the pre-last step we use the inclusion $F_A \in \mathcal{R}_{-\alpha}$, due to $F \in MRV(\alpha,\, V,\,\mu)$, $\mu(A)>0$, see Remark \ref{rem.KP.3.A}. Using representation of \eqref{eq.KP.2.11} on \eqref{eq.KP.4.14}, we find \eqref{eq.KP.4.7}.~\halmos
\end{enumerate}

Next, we will establish the asymptotic behavior of finite-time ruin probability in the model described by \eqref{eq.KP.4.18}. In order to define the ruin probability in this model, we need a family of ruin-sets, that are given in \cite[Ass. 5.1]{samorodnitsky:sun:2016}. We say that a set $\bbb$ is 'decreasing' if the $-\bbb$ is 'increasing' and by $\partial \bbb$ we denote the border of set $\bbb$.

\begin{assumption} \label{ass.KP.4.2}
Let $L$ be a ruin-set, which is an open, decreasing set, such that ${\bf 0} \in \partial L$, the complement set $L^c$ is convex, and $L=x\,L$, for any $x>0$.
\end{assumption}

According to this assumption, the finite time ruin probability is defined as follows:
\beam \label{eq.KP.4.19}
\psi_{{\bf l},L}(x,\,t) = \PP\left[ {\bf U}(s) \in L\,, \;\exists \; s\in [0,\,t] \right]\,,
\eeam
From \eqref{eq.KP.4.19} we can see that the finite time ruin probability, in addition to the initial capital, is also depends on the ruin-set and the allocation of the capital through the weights $(l_1,\,\ldots,\,l_d)=:{\bf l}$. Two important forms of ruin-set, satisfying Assumption  \ref{ass.KP.4.2}, are 
\beao
L_1=\{{\bf x}\;:\; x_i <0\,, \;\text{ for some} \; i=1,\,\ldots,\,d \}\,, \quad L_2=\left\{{\bf x}\;:\; \sum_{i=1}^d x_i <0 \right\}\,,
\eeao
where in $L_1$ is expressed the ruin, is due to negative surplus in one line of business, while in $L_2$ the ruin is due to negative total surplus. In multivariate risk models there are several kinds of ruin probabilities, see for example \cite{cheng:yu:2019}. Let us denote
\beao
\int_0^te^{-r\,y}\,{\bf c}(y)\,dy:= \left(
\begin{array}{c}
\int_{0}^t e^{-r\,y}\,c_1(y)\,dy \\ 
\vdots \\ 
\int_{0}^t e^{-r\,y}\,c_d(y)\,dy 
\end{array} 
\right)\,,
\eeao
Then, the set $A={\bf l}-L$ belongs to $\mathscr{R}$, see \cite[Sec. 5]{samorodnitsky:sun:2016}. Therefore, relation \eqref{eq.KP.4.19} is equivalent to  
\beao
\psi_{{\bf l},L}(x,\,t) &=& \PP\left[ {\bf U}(s) \in L\,, \;\exists \; s\in [0,\,t] \right]\\[2mm]
&&= \PP\left[ {\bf D}_r(s) -\int_0^s e^{-r\,y}\,{\bf c}(y)\,dy \in x\,({\bf l} - L)\,,\;\exists\; s\in [0,\,t] \right]\\[2mm] 
&&= \PP\left[ {\bf D}_r(s) -\int_0^s e^{-r\,y}\,{\bf c}(y)\,dy \in x\,A\,,\;\exists\; s\in [0,\,t] \right]\,,
\eeao
where in the second step we used the relation $x\,L = L$, from Assumption \ref{ass.KP.4.2}.

\bco \label{cor.KP.4.2}
Let some fixed ruin-set $L$, that satisfies Assumption \eqref{ass.KP.4.2} and $A={\bf l} - L$. In risk model \eqref{eq.KP.4.18}
\begin{enumerate}
\item[(i)]
If for any fixed $t \in (0,\,\infty)$ the conditions of Theorem \ref{th.KP.4.1}(i) are true, then it holds
\beam \label{eq.KP.4.23}
\psi_{{\bf l},L}(x,\,t) \sim C_0\, \PP\left[ {\bf X}\in x\,A \right]\,.
\eeam
\item[(ii)]
If for any fixed $t \in (0,\,\infty)$ the conditions of Theorem \ref{th.KP.4.1}(ii) are true, then it holds
\beam \label{eq.KP.4.26}
\psi_{{\bf l},L}(x,\,t) \sim C_r\,\,\mu(A)\,\overline{V}(x)\,.
\eeam
\end{enumerate}
\eco

\pr~
Let consider $r\geq 0$. By Theorem \ref{th.KP.4.1} we obtain relations \eqref{eq.KP.4.5} and \eqref{eq.KP.4.7} , with $C_0,\,C_r \in (0,\,\infty)$. Hence, by closure property of class $\mathcal{S}_A$ with respect to strong equivalence, see \cite[Prop. 4.12(a)]{samorodnitsky:sun:2016}, we find out that the ${\bf D}_r(s)$ belongs to class $\mathcal{S}_A$ in both cases (see also \eqref{eq.KP.2.19}). Hence, because $xA$ is increasing set, and ${\bf D}_r(s)$ is non-decreasing with respect to time, we have
\beam \label{eq.KP.4.24} \notag
\psi_{{\bf l},L}(x,\,t) &\leq& \PP\left[ {\bf D}_r(t) \in x\,A + \int_0^s e^{-r\,y}\,{\bf c}(y)\,dy \,,\;\exists\; s\in [0,\,t] \right]\\[2mm]
&\leq & \PP\left[ {\bf D}_r(t) \in x\,A \right] \sim C_r\, \PP\left[ {\bf X}\in x\,A \right]\,.
\eeam
On the other hand side, because of upper bound of the premium density, and from the subexponential behavior of ${\bf D}_r(s)$, through the \cite[Prop. 4.12(b)]{samorodnitsky:sun:2016}, we get
\beam \label{eq.KP.4.25} \notag
\psi_{{\bf l},L}(x,\,t) &\geq& \PP\left[ {\bf D}_r(t) - \int_0^t e^{-r\,y}\,{\bf c}(y)\,dy  \in x\,A \right]\\[2mm]
&\geq & \PP\left[ {\bf D}_r(t) \in x\,A + {\bf M}\,t \right] \sim \PP\left[ {\bf D}_r(t) \in x\,A \right] \sim C_r\, \PP\left[ {\bf X}\in x\,A \right]\,,
\eeam
where ${\bf M}=(M_1,\,\ldots,\,M_d)^{\top}$. Therefore by relations \eqref{eq.KP.4.24} and \eqref{eq.KP.4.25} we reach the desired \eqref{eq.KP.4.26} (and relation \eqref{eq.KP.4.23} in the first case, when $r=0$).
~\halmos

\section{Vector type precise large deviations} \label{sec.KP.7}

The precise large deviations in non-random and random sums of heavy tailed random variables, are well-studied, see for example in \cite{mikosch:nagaev:1998}, \cite{ng:tang:yan:yang:2004}, \cite{konstantinides:mikosch:2005}, \cite{wang:wang:cheng:2006} among many others. The known problem of lower bound of precise large deviation can be expressed as
\beam \label{eq.KP.5.1} 
\liminf_{\nto} \inf_{x\geq \gamma\,n} \dfrac{\PP\left[ \sum_{i=1}^n Z_i > x +n\,c \right]}{\sum_{i=1}^n \overline{V}_i(x+n\,c)} \geq 1 \,,
\eeam
for any $\gamma > \gamma_0$, where $\gamma_0 =\gamma_0(c)$ and $c$ some appropriately chosen  constant, where $Z_i \stackrel{d}{\sim} V_i$, with $i \in \bbn$. Relation \eqref{eq.KP.5.1} was examined in several papers, see for example in \cite{konstantinides:loukissas:2011}, \cite{loukissas:2012} for identically distributed $Z_i \stackrel{d}{\sim} V \in \mathcal{L}$ and weakly dependent or independent cases respectively. In \cite{he:cheng:wang:2013} relation \eqref{eq.KP.5.1} was established for non-identically distributed summands, under a more general dependence structure, and with only condition on summand distributions the infinite right endpoint. In multivariate set up the precise large deviations are studied under the following concept. For some integer $d \in \bbn$ the sequences $\{X_{1j}\,,\;j\geq 1\},\,\ldots,\,\{X_{dj}\,,\;j\geq 1\}$ are mutually independent and the focus was put on the study of the probability
\beam \label{eq.KP.5.2} 
\PP\left[ S_d > x\right] \,,
\eeam
as $\nto$, or similar (centralized) quantities, where we denote
\beao
S_d = \sum_{i=1}^d S_{n_i}\,, \qquad \qquad S_{n_i}= \sum_{j=1}^{n_i} X_{ij}\,,
\eeao
with $n_i \in \bbn$, for any $i=1,\,\ldots,\,d$, see for example \cite{wang:wang:2007}, \cite{lu:2012}, \cite{chen:cheng:2024b}. Although this approach of multivariate models contains several positive features, as for example that the number of summands $n_i$, $n_j$ is NOT necessarily the same, for $1\leq i \neq j \leq d$, and that the distributions belong to rather large classes, like $\mathcal{D}$ or $\mathcal{L}$, there are two gaps. The first gap is, that relation \eqref{eq.KP.5.2} focus only to the total sum of the $d$-sequences and it does not bring any information for each sequence separately, namely it is only connected with the sum-type ruin probability. The second one is, that is required independence among the sequences and each sequence permits only weakly dependent terms, with obvious restrictions in dependence modeling. From the other side in the frame of multivariate random walks, see \cite{hult:lindskog:mikosch:samorodnitsky:2005}, the focus was in large deviations on some rare-sets, when the random walk steps are independent and identically distributed random vectors with $MRV$ distribution, where each vector has arbitrarily dependent components. In \cite{haegele:lehtomaa:2021} was examined similar large deviations, in a distribution class  that does not belong to $MRV$. This second approach have also some difficulties to capture the dependence between the vectors, as also the fact that the large deviations are not such practical as the precise ones.

\subsection{Lower bounds of precise large deviations} \label{sec.KP.7.a}

In our attempt to cover the gaps of these two previous approaches, we establish precise asymptotic bounds for the large deviations on the sets $A$ from the family of sets $\mathscr{R}$. The first result estimates the lower bound of precise large deviations of non-random sums
\beao
{\bf S}_n = \sum_{i=1}^n {\bf X}^{(i)}\,,
\eeao
for any $n \in \bbn$, with ${\bf S}_0={\bf 0}$. The distributions of the random vectors ${\bf X}^{(i)}$ should have marginals with infinite right endpoint (more precisely, only $F_A^{(i)}$ should have infinite right endpoint). Furthermore, each vector has arbitrarily dependent components and the random vectors are mutually weakly dependent, containing a dependence structure stemming from \cite{he:cheng:wang:2013}. We should mention that the next result coincides to \cite[Theorem 2.1]{he:cheng:wang:2013}, when $d=1$ and $A=(1,\,\infty)$. 

\bth \label{th.KP.5.1}
Let $A \in \mathscr{R}$ be some fixed set. We consider the sequence $\{ {\bf X}^{(i)} \,,\; i\in \bbn \}$ of non-negative random vectors, with distributions 
$\{ F_{i}\,,\; i\in \bbn \}$ respectively. We assume that it holds
\beam \label{eq.KP.5.5} 
\lim_{\nto} \sup_{x \geq \gamma\,n} \sup_{1\leq i < j \leq n} x\,\PP\left[{\bf X}^{(i)} \in x\,A\;|\; {\bf X}^{(j)} \in x\,A \right] =0\,,
\eeam
for any $\gamma >0$. Then we obtain
\beam \label{eq.KP.5.6} 
\liminf_{\nto} \inf_{x \geq \gamma\,n} \dfrac{\PP\left[{\bf S}_n \in x\,A \right]}{\sum_{i=1}^n \PP\left[{\bf X}^{(i)} \in x\,A \right] } \geq 1\,,
\eeam
for any $\gamma >0$.
\ethe

\pr~
Since the ${\bf X}^{(i)}$ are non-negative random vectors and the set $xA$ is increasing, it follows that for any $x>0$ it holds
\beam \label{eq.KP.5.7}
\PP\left[{\bf S}_n \in x\,A \right] &\geq& \PP\left[\bigcup_{i=1}^n \{{\bf X}^{(i)} \in x\,A \}\right] \\[2mm] \notag
&\geq& \sum_{i=1}^{n} \PP\left[{\bf X}^{(i)} \in x\,A \right]- \sum_{1\leq i < j \leq n} \PP\left[{\bf X}^{(i)} \in x\,A\,,\; {\bf X}^{(j)} \in x\,A\right] \,.
\eeam
For the second sum in the right hand side of relation \eqref{eq.KP.5.7}, from relation \eqref{eq.KP.5.5}, we find that for any arbitrarily chosen $\vep \in (0,\,1)$, there exists some $n_0=n_0(\vep, \gamma)>0$, such that for any $n\geq n_0$, it holds
\beam \label{eq.KP.5.8} \notag
&&\sum_{1\leq i < j \leq n} \PP\left[{\bf X}^{(i)} \in x\,A\,,\; {\bf X}^{(j)} \in x\,A\right]= \sum_{j=2}^n \PP\left[{\bf X}^{(j)} \in x\,A\right] \sum_{i=1}^{j-1} \PP\left[{\bf X}^{(i)} \in x\,A\;|\; {\bf X}^{(j)} \in x\,A\right] \\[2mm] 
&&<\vep\,\gamma\,\sum_{j=2}^n \PP\left[{\bf X}^{(j)} \in x\,A\right] \,\dfrac {n-1}x<\vep\,\sum_{j=2}^n \PP\left[{\bf X}^{(j)} \in x\,A\right] <\vep\, \sum_{i=1}^{n} \PP\left[{\bf X}^{(i)} \in x\,A \right]\,,
\eeam
uniformly, for any $x\geq \gamma\,n$, where in the pre-last step we use the inequality $\gamma\,(n-1)<x$. By relations \eqref{eq.KP.5.7} and \eqref{eq.KP.5.8}, leaving the $\vep$ to tend to zero, we find \eqref{eq.KP.5.6}. 
~\halmos

\bre \label{rem.KP.5.1}
Let $c$ be a real number. Then we can verify that relation \eqref{eq.KP.5.6}, is equivalent to 
\beam \label{eq.KP.6.9}
\liminf_{\nto} \inf_{x \geq \gamma\,n} \dfrac{\PP\left[{\bf S}_n \in (x+n\,c)\,A \right]}{\sum_{i=1}^n \PP\left[{\bf X}^{(i)} \in (x+n\,c)\,A \right] } \geq 1\,,
\eeam 
for $\gamma > -c$. Relation \eqref{eq.KP.6.9} plays crucial role in risk theory, since by $c$ is understood the premium rate in insurance business, while by ${\bf X}^{(1)},\,\ldots,\,{\bf X}^{(n)}$ are understood the multivariate claims. In the case of set $A$ of the form 
$A=\left\{ {\bf y}\;:\;\sum_{i=1}^d k_i\,y_i > 1 \right\}$, with $k_i\geq 0$ for any $i=1,\,\ldots,\,d$ and $k_1 +\cdots + k_d=1$. In such a case $x$ is understood as the initial capital, and $c$ is the sum of premium rates of the $d$-lines of business. Then relation  \eqref{eq.KP.6.9} is connected with the sum-type ruin probability.
\ere

Let us focus on the lower bound of the precise large deviations for the random sum ${\bf S}_{N(t)} = \sum_{i=1}^{N(t)} {\bf X}^{(i)}$, with $\{N(t)\,,\;t\geq 0 \}$ ($N(0)=0$) is a counting process, independent of $\{{\bf X}^{(i)}\,,\;i \in \bbn \}$. Let us denote by $\lambda(t) = \E[N(t)] < \infty$, for any fixed $t\geq 0$, with $\lambda (t) \to \infty$, as $t\to \infty$, its renewal mean. We suppose that the counting process $\{N(t)\,, \;t\geq 0\}$, satisfies the condition
\beam \label{eq.KP.6.11}
\dfrac {N(t)}{\lambda(t) } \stackrel{P}{\rightarrow} 1\,,
\eeam
as $t\to \infty$, that represents a standard condition in one-dimensional precise large deviations, see for example \cite{kluppelberg:mikosch:1997}, 
\cite{konstantinides:loukissas:2011}, \cite{loukissas:2012} among others. According to comments from \cite{kluppelberg:mikosch:1997}, relation \eqref{eq.KP.6.11} is satisfied by all the renewal processes, as well as by the inhomogeneous Poisson processes, therefore represents a rather general condition.

Inspired by \cite[Th. 3.1]{he:cheng:wang:2013}, the next theorem  contain it as special case, when $d=1$ and $A=(1,\,\infty)$. In both cases, the random vectors  $\{{\bf X}^{(i)}\,,\;i \in \bbn \}$ should satisfy some kind of weak-equivalence on sets $xA$. We formulate this restriction in the following assumption.

\begin{assumption} \label{ass.KP.6.1}
Let $A \in \mathscr{R}$ be some fixed set. We consider a non-negative random vector ${\bf X}$, such that it holds
\beao 
c_1 =\liminf \inf_{i\geq 1} \dfrac{\PP\left[{\bf X}^{(i)}\in x\,A\right]}{\PP\left[{\bf X}\in x\,A\right]}\leq \limsup \sup_{i\geq 1} \dfrac{\PP\left[{\bf X}^{(i)}\in x\,A\right]}{\PP\left[{\bf X}\in x\,A\right]} =c_2\,,
\eeao
with $0<c_1 \leq c_2 < \infty$.
\end{assumption}

\bth \label{th.KP.6.2}
Let $A \in \mathscr{R}$ be some fixed set. We consider the sequence  $\{{\bf X}^{(i)}\,,\;i \in \bbn \}$ of non-negative random vectors, with distributions  
$\{F_{i}\,,\;i \in \bbn \}$, respectively. We suppose that the  $\{{\bf X}^{(i)}\,,\;i \in \bbn \}$ satisfy condition  \eqref{eq.KP.5.5} and Assumption \ref{ass.KP.6.1}. We assume also that the $\{N(t)\,,\;t\geq 0 \}$ is a counting process, independent of  $\{{\bf X}^{(i)}\,,\;i \in \bbn \}$, that satisfies 
relation \eqref{eq.KP.6.11}. Then it holds
\beam \label{eq.KP.6.13}
\liminf_{\tto} \inf_{x \geq \gamma\,\lambda(t)} \dfrac{\PP\left[{\bf S}_{N(t)}\in x\,A \right]}{\sum_{i=1}^{\left\lfloor \lambda(t) \right\rfloor} \PP\left[{\bf X}^{(i)} \in x\,A \right] } \geq 1\,,
\eeam 
for any  $\gamma >0$.
\ethe

\pr~
At first, we see that when $x>0$, then $\PP\left[{\bf S}_0 \in x\,A\right]=0$. Further, by relation  \eqref{eq.KP.6.11} for the counting process, for any $\delta,\,\vep \in (0,\,1)$ and for any $\gamma >0$, there exists some $t'=t'(\delta,\,\vep)>0$, such that when $t\geq t'$, it holds uniformly for $x \geq \gamma\,\lambda(t)$ the following inequality
\beam \label{eq.KP.6.14} \notag
&&\PP\left[{\bf S}_{N(t)} \in x\,A \right]=\sum_{n=1}^{\infty} \PP\left[{\bf S}_{n} \in x\,A \right]\,\PP[N(t)=n] \geq \sum_{n=\left\lfloor (1-\delta)\lambda(t) \right\rfloor}^{\left\lfloor(1+\delta) \lambda(t) \right\rfloor} \PP\left[{\bf S}_{n} \in x\,A \right]\,\PP[N(t)=n] \\[2mm]
&&\geq \PP\left[{\bf S}_{\left\lfloor (1-\delta)\lambda(t) \right\rfloor} \in x\,A \right]\, \sum_{n=\left\lfloor (1-\delta)\lambda(t) \right\rfloor}^{\left\lfloor(1+\delta) \lambda(t) \right\rfloor}\PP[N(t)=n] \\[2mm] \notag
&&\geq \PP\left[{\bf S}_{\left\lfloor (1-\delta)\lambda(t) \right\rfloor} \in x\,A \right]\,\PP\left[ \left|\dfrac{N(t)}{\lambda(t)} -1\right|< \delta \right]\geq (1-\vep)\,\PP\left[{\bf S}_{\left\lfloor (1-\delta)\lambda(t) \right\rfloor} \in x\,A \right]\,,
\eeam
where in the third step, we used that the $\left\{{\bf X}^{(i)}\,,\;i \in \bbn \right\}$ are non-negative random vectors, and the set $x\,A$ is increasing, see \eqref{eq.KP.2.14}, \eqref{eq.KP.2.15}.

By relation \eqref{eq.KP.6.14}, taking into account Theorem \ref{th.KP.5.1}, we can find some $\widetilde{t}= \widetilde{t}(\delta,\,\vep,\,\gamma) \geq t'$, such that when $t>\widetilde{t}$, it holds uniformly for $x \geq \gamma\,\lambda(t)$ that,
\beam \label{eq.KP.6.15} 
\PP\left[{\bf S}_{N(t)} \in x\,A \right] \geq (1-\vep)^2\,\sum_{i=1}^{\left\lfloor (1-\delta)\lambda(t) \right\rfloor} \PP\left[{\bf X}^{(i)} \in x\,A \right]\,,
\eeam
 By Assumption \ref{ass.KP.6.1}, we find that there exists some $t^{*}=t^{*}(\delta,\,\vep,\,\gamma)\geq \widetilde{t}$, such that when $t>t^{*}$, it holds uniformly for $x \geq \gamma\,\lambda(t)$ that,
\beam \label{eq.KP.6.16} \notag
\sum_{i=1}^{\left\lfloor (1-\delta)\lambda(t) \right\rfloor} \PP\left[{\bf X}^{(i)} \in x A \right]&=&\sum_{i=1}^{\left\lfloor \lambda(t) \right\rfloor} \PP\left[{\bf X}^{(i)} \in x A \right]\,\left(1-\dfrac{\sum_{i={\left\lfloor (1-\delta)\lambda(t) \right\rfloor +1}}^{\left\lfloor \lambda(t) \right\rfloor} \PP\left[{\bf X}^{(i)} \in x A \right]}{\sum_{i=1}^{\left\lfloor \lambda(t) \right\rfloor} \PP\left[{\bf X}^{(i)} \in x A \right]} \right)\\[2mm] 
&\geq& \left(1-2\,\delta\,\dfrac{c_2}{c_1} \right)\,\sum_{i=1}^{\left\lfloor \lambda(t) \right\rfloor} \PP\left[{\bf X}^{(i)} \in x\,A \right]\,,
\eeam
where $c_1,c_2$ are defined in Assumption \eqref{ass.KP.6.1}. From relations \eqref{eq.KP.6.15}, \eqref{eq.KP.6.16} and the arbitrary choice of $\vep$, $\delta$, we conclude \eqref{eq.KP.6.13}.
~\halmos
 
\bre \label{rem.KP.6.2}
As in Remark \ref{rem.KP.5.1},similarly in the case of random sums, if $c$ be a real number, then relation \eqref{eq.KP.6.13} is equivalent to 
\beao
\liminf_{t \to \infty} \inf_{x \geq \gamma\,\lambda(t)} \dfrac{\PP\left[{\bf S}_{N(t)} \in (x+c\,\lambda(t))\,A \right]}{\sum_{i=1}^{\left\lfloor \lambda(t) \right\rfloor} \PP\left[{\bf X}^{(i)} \in (x+c\,\lambda(t))\,A \right] } \geq 1\,,
\eeao
for $\gamma > -c$.
\ere

\subsection{An application to compound risk model} \label{sec.KP.7.b}

Here, we provide an application of Theorem \ref{th.KP.6.2} with respect to lower bound of precise large deviations in a compound risk model. We consider that the accidents arrive at the time moments $\{ \tau_i\,,\; i \in \bbn\}$, with $\tau_0=0$, which represent a counting process 
\beam \label{eq.KP.6.17}
\tau(t) = \sup\{ n\in \bbn\;:\; \tau_n \leq t\}\,,
\eeam 
with mean value $\nu(t) := \E[\tau(t)]$. In several applications, as for example in risk theory, is natural to assume that (in each portfolio) are produced more than one claims at each accident arrival $\tau_i$, with $i \in \bbn$. Hence, we assume that at each $\{ \tau_i\,,\; i \in \bbn\}$ produced $\{ Z_i\,,\; i \in \bbn\}$ number of claim-vectors, that are integer valued random variables, non-degenerate at zero. Then for some time $t \geq 0$, the number of claim vectors is given through the process
\beam \label{eq.KP.6.18}
M(t) = \sum_{i=1}^{\tau(t)} Z_i\,,
\eeam 
while the aggregate claims up to this time  $t \geq 0$, are provided by the relation
\beam \label{eq.KP.6.19}
{\bf S}_{M(t)} = \sum_{i=1}^{M(t)} {\bf X}^{(i)}=\left(\sum_{i=1}^{M(t)} X_{1}^{(i)}\,,\ldots,\,\sum_{i=1}^{M(t)} X_{d}^{(i)} \right)^{\top}\,.
\eeam
As in Section 5, each ${\bf X}^{(i)}$ can have zero components, but not all of the components equal to zero. In one-dimensional case, the risk model \eqref{eq.KP.6.19} was introduced by \cite{tang:su:jiang:zhang:2001}, and was studied later by many papers, see \cite{aleskeviciene:leipus:siaulys:2008}, \cite{konstantinides:loukissas:2010}, \cite{loukissas:2012}, \cite{zhang:yang:xu:2026}, among others. In relation to precise large deviations we find that in \cite{konstantinides:loukissas:2010}, was generalized the result by \cite{tang:su:jiang:zhang:2001}, while in \cite{konstantinides:loukissas:2010} was given the lower bound for the class $\mathcal{L}$. Although the last papers consider the Assumption \ref{ass.KPX.6.3} (for one-dimensional claims), they considered that the $\{\tau(t)\,,\; t \geq 0\}$ represents a renewal process, which we weaken here, see Assumption  \ref{ass.KPX.6.4} below.
	
\begin{assumption} \label{ass.KPX.6.3}
Let $\{Z_{i}\,,\; i \in \bbn \}$ be independent, identically distributed random variables, with mean $\E[Z_1]< \infty$, while $\{\Theta_{i}:=\tau_{i} - \tau_{i-1}\,,\; i \in \bbn \}$ be identically distributed random variables, with mean $\E[\Theta_1]=1/\lambda \in (0,\,\infty)$. Further we also assume that
$\{Z_{i}\,,\; i \in \bbn \}$, $\{\Theta_{i}\,,\; i \in \bbn \}$ and $\{{\bf X}^{(i)}\,,\;i \in \bbn \}$ are mutually independent. 
\end{assumption}

In the following assumption we consider that the dependence structure among the interarrival times  $\{\Theta_{i}\,,\; i \in \bbn \}$, is the extended negatively orthant dependence, symbolically $ENOD$. We refer to \cite{wang:cheng:2011} and \cite{wang:cui:wang:ma:2012} for more discussions and applications of $ENOD$.

\begin{assumption} \label{ass.KPX.6.4}
Suppose that for the $\{\Theta_{i}\,,\; i \in \bbn \}$, there exists some constant $M>0$, such that for any $x_i \in \bbr$, and any $i=1,\,\ldots,\,n$, with $n \in \bbn$, hold the inequalities
\beam \label{eq.KP.6.20}
\PP \left( \bigcap_{i=1}^n \{\Theta_i > x_i \} \right) \leq M\,\prod_{i=1}^n \PP \left(\Theta_i > x_i \right)\,\,, 
\PP \left( \bigcap_{i=1}^n \{\Theta_i \leq x_i \} \right) \leq M\,\prod_{i=1}^n \PP \left(\Theta_i \leq x_i \right)\,,
\eeam
\end{assumption}

The result of the following corollary remains new even for the one-dimensional subcase with $A=(1,\,\infty)$, since it generalizes \cite[Prop. 4.2]{loukissas:2012} (in the centralized form).

\bco \label{cor.KPX.6.1}
Let $A \in \mathscr{R}$ be some fixed set and let the compound risk model given by \eqref{eq.KP.6.19}, satisfy Assumptions \ref{ass.KPX.6.3} and \ref{ass.KPX.6.4}. Moreover, let $\{{\bf X}^{(i)}\,,\;i \in \bbn \}$ has distributions  $\{F_{i}\,,\; i \in \bbn \}$ satisfying the dependence structure \eqref{eq.KP.5.5} and Assumption \ref{ass.KP.6.1}. Then for any $\gamma>0$ it holds
\beam \label{eq.KP.6.21}
\liminf_{t\to \infty} \inf_{x \geq \gamma \lambda t \E[Z]} \dfrac{\PP \left( {\bf S}_{M(t)} \in x\,A \right)}{\sum_{i=1}^{\left\lfloor \lambda t \E[Z]\right\rfloor}\PP \left( {\bf X}^{(i)} \in x\,A \right)} \geq 1\,.
\eeam
\eco

\pr~
Firstly, from Assumption \ref{ass.KPX.6.4}, through \cite[Th. 1.4]{wang:cheng:2011} we obtain $\tau(t) \sim \lambda\,t$, as $t \to \infty$, almost surely, while by \cite[Th. 1.2]{wang:cheng:2011} we find $\nu(t)=\E[\tau(t)] \sim \lambda\,t$, as $t \to \infty$, hence relation \eqref{eq.KP.6.11} is satisfied. Thus, by Assumption \ref{ass.KPX.6.3} it holds: 
\beao
\E[M(t)] = \E \left[\sum_{i=1}^{\tau(t)} Z_i \right]= \nu(t)\,\E[Z] \sim \lambda\,t\,\E[Z]\,,
\eeao 
as $\tto$. Thus, relation \eqref{eq.KP.6.21} is directly implied by application of Theorem \ref{th.KP.6.2} on 
\eqref{eq.KP.6.19}. 
~\halmos

\noindent \textbf{Acknowledgments.} 
We would like to thank Prof. Yuebao Wang, for his useful advice that improved substantially the text. We are also very grateful to the anonymous referees, for their constructive comments on an earlier version of the paper.

\end{document}